\documentclass[11pt]{article}

\usepackage{latexsym,epsfig}
\usepackage{color}
\usepackage{amssymb,amsmath,amsfonts}
\usepackage{exscale}
\usepackage{ifthen}
\usepackage{longtable}
\usepackage{subfigure}
\usepackage{caption}
\usepackage{cite}
\usepackage{lscape}
\usepackage{enumitem}
\newcommand{\vol}{\mbox{\rm vol}}

\newcommand{\R}{\mathbb{R}}

\newtheorem{definition}{\bf Definition}[section]
\newtheorem{definitions}[definition]{\bf Definitions}
\newtheorem{corollary}[definition]{\bf Corollary}
\newtheorem{theorem}[definition]{\bf Theorem}
\newtheorem{lemma}[definition]{\bf Lemma}
\newtheorem{proposition}[definition]{\bf Proposition}
\newtheorem{example}[definition]{\bf Example}
\newtheorem{examples}[definition]{\bf Examples}
\newtheorem{remark}[definition]{\sc Remark}
\newtheorem{applemma}[definition]{\bf Approximation Lemma}

\begin{document}
\bibliographystyle{alphabetical}
\protect\pagenumbering{arabic}

\title{\sc Asymptotic geometry  of negatively curved manifolds of finite volume}
\author{\sc F. Dal'Bo, M. Peign\'e,
 J.C. Picaud \& A. Sambusetti}
\date{\today}
\maketitle

\centerline {\small {\bf Abstract}}
We study the asymptotic behaviour of simply connected, Riemannian manifolds $X$ of strictly negative curvature admitting a non-uniform lattice $\Gamma$. 
If the quotient manifold $\bar X= \Gamma \backslash X$ is  asymptotically $1/4$-pinched, we prove that  $\Gamma$  is divergent and $U\bar X$ has finite Bowen-Margulis measure (which is then ergodic and totally conservative with respect to the geodesic flow); moreover, we show  that, in this case,  the volume growth of balls $B(x,R)$ in $X$ is asymptotically equivalent to a purely exponential function $c(x)e^{\delta R}$, where $\delta$ is the topological entropy of the geodesic flow of $\bar X$. \linebreak This generalizes Margulis' celebrated theorem to negatively curved spaces of finite volume. 
In contrast,  we exhibit examples  of lattices $\Gamma$ in negatively curved spaces $X$ (not asymptotically $1/4$-pinched)    where, depending on the critical exponent of the parabolic subgroups   and on the finiteness of the Bowen-Margulis measure, the growth function is exponential, lower-exponential or even upper-exponential. 
\\


 \section{Introduction}
Let $X$  be a  complete, simply connected manifold with strictly negative curvature. \linebreak In the sixties,  G. Margulis \cite{margulis},  using measure theory on the foliations of the Anosov system defined by the geodesic flow, showed that if $\Gamma$ is   a uniform lattice of $X$ (i.e.  a torsionless,  discrete group of isometries  such that $\bar X = \Gamma \backslash X$ is compact), then the  {\em orbital function} of $\Gamma$ is asymptotically equivalent to a purely exponential function:
$$v_{\Gamma}(x,y,R) = \#  \{ \gamma \in \Gamma \; | \; d(x, \gamma y) < R \} \;\sim \; c_\Gamma (x,y) e^{\delta(\Gamma) R}$$
where $\delta(\Gamma)=\lim_{R \rightarrow \infty}  R^{-1} \ln v_\gamma (x,x,R)$ is  the {\em critical exponent} of $\Gamma$, and $\sim$ means that the quotient tends to $1$ when $R \rightarrow \infty$.
By integration over fundamental domains, one then obtains  an asymptotic equivalence for the {\em volume growth function} of $X$:  
$$v_X(x, R) = \vol B(x,R) \;\sim \; m(x) e^{\delta(\Gamma)R} \,.$$

\noindent It is well-known that the exponent  $\delta(\Gamma)$  equals the {\em topological entropy} of the geodesic flow of $\bar X$  (see \cite{otalpeigne}) and that, for uniform lattices,  it is the same as  the {\em volume entropy} $\omega (X) = \limsup \frac{1}{R} \ln v_X(x, R)$ of the manifold $X$.   The function $m(x)$, depending on the center   of the ball, is the {\em Margulis function} of $X$.

Since then, this result has been generalized in different directions. Notably, G. Knieper showed in \cite{knieper} that the volume growth function of a Hadamard space $X$  (a complete, simply connected manifolds with nonpositive curvature) admitting uniform lattices is {\em purely exponential}, provided  that $X$ has rank one, that is $v_X(x, R) \asymp   e^{ \omega (X)R} $ 
 (where $f \asymp  g $ means that $1/A < f(R)/g(R) < A$ for some positive $A$, when $R\gg0$). \linebreak 
In general, he showed that  $v_X(x, R) \asymp   R^{\frac{d-1}{2}}e^{\omega (X) R}$ for rank $d$ manifolds; however, as far as the authors are aware, it is still unknown whether there exists a Margulis function for Hadamard manifolds of rank 1 with uniform lattices, i.e. a function $m(x)$ such that  $v_X(x, R) \sim m(x) e^{ \omega (X) R}$, even in the case of surfaces. 
Another remarkable case is that of {\em asymptotically harmonic manifolds} of strictly negative curvature, where the strong asymptotic homogeneity implies the existence of a Margulis function, even without compact quotients, cp.\cite{CS}.

In another direction, it seems   natural   to ask what happens for a Hadamard space  $X$ of negative curvature admitting  {\em nonuniform lattices} $\Gamma$  (i.e. $\vol (\Gamma \backslash X) < \infty$): {\em is  $v_X$  purely exponential and, more precisely, does $X$  admit a Margulis function?} \linebreak 
Let us  emphasize  that if $X$ also admits a uniform lattice  then $X$ is a symmetric space  of rank one (by \cite{eberlein}, Corollary 9.2.2); therefore, we are  interested in spaces  which do not have uniform lattices, i.e. the universal covering  of finite volume, negatively curved  manifolds which are not locally symmetric.

 It is worth to stress here that 
the orbital function of   $\Gamma$   is closely related to  the volume growth function of $X$, but it generally has, even for lattices,  a different asymptotic behaviour than $v_X(x,R)$.  
A precise asymptotic equivalence fo $v_{\Gamma}$ was proved by  T. Roblin \cite{roblin}  in a very general setting.  Namely, he proved that for any  nonelementary group of isometries $\Gamma$ of a CAT(-1) space $X$ with non-arithmetic  length spectrum\footnote{This means that the additive subgroup of $\R$  generated by the length of closed geodesics in $\bar X= \Gamma \backslash X$   is dense in $\R$; it is the case, for instance, if $dim(X)=2$, or when $\Gamma$ is a lattice.} and $\bar X = \Gamma \backslash X$, one has:

{\em 
\noindent (a)  $v_{\Gamma}(x,y,R) \sim c_\Gamma (x,y)  e^{\delta (\Gamma)R}$ if the Bowen-Margulis measure on the unitary tangent bundle $U\bar X$ is finite;

\noindent (b) $v_{\Gamma}(x,y,R) = o(R) e^{\delta (\Gamma)R}$, where $o(R)$ is infinitesimal, otherwise.
}

\noindent Thus, the behaviour of $v_\Gamma (x,R)$ strongly depends on the finiteness of the Bowen-Margulis measure $\mu_{BM}$; also,  the asymptotic constant can be expressed  in terms of $\mu_{BM}$ and of the family of Patterson-Sullivan  measures $(\mu_x)$ of $\Gamma$, as $c_\Gamma (x,y) = \frac{\parallel \! \mu_x \! \parallel \, \parallel \! \mu_y \! \parallel}{ \delta (\Gamma) \cdot \parallel \! \mu_{BM} \! \parallel}$. \linebreak
 In section \S4 we will recall a useful criterion ({\em Finiteness Criterion} (\ref{eqcritfiniteness}), due to Dal'Bo-Otal-Peign\'e), ensuring  that a geometrically finite group has $\mu_{BM} (U\bar X) < \infty$, hence a precise asymptotics for $v_{\Gamma}(x,R)$ as in (a).
    
    \pagebreak
         
 \noindent On the other hand, any convergent group $\Gamma$ exhibits  a behaviour as in (b), since it certainly has  infinite Bowen-Margulis measure (by Poincar\'e recurrence,  $\mu_{BM} (U\bar X)<\infty$ \linebreak implies that the geodesic flow is totally conservative, and this is equivalent to divergence, by Hopf-Tsuji-Sullivan's theorem).  
Notice that, whereas uniform lattices always are divergent and with finite Bowen-Margulis measure,   for nonuniform lattices  $\Gamma$  divergence  and condition (\ref{eqcritfiniteness}) in general may fail. Namely, this can happen  only in case  $\Gamma$ has a ``very  large'' parabolic subgroup $P$,  that is such that $\delta(P) = \delta(\Gamma)$:  we   will call  {\em exotic} such  a lattice $\Gamma$, and we will say that such a  $P$ is a {\em dominant} parabolic subgroup. \linebreak
 Convergent, exotic lattices are constructed by the authors in \cite{DPPSnew}; also, one can  find in \cite{DPPSnew} some original counting results for the orbital function of $\Gamma$ in infinite Bowen-Margulis measure, more precise than (b).

 However, as we shall see, the volume growth function $v_X$  has  a wilder behaviour than \nolinebreak $v_\Gamma$.
 In \cite{DPPS}  we  proved  that    for nonuniform  lattices in pinched, negatively curved spaces $X$, the  functions $v_{\Gamma}$ and $v_X$ can have different exponential growth rates, i.e. $\omega (X) \neq \delta(\Gamma)$. In the Example \ref{exsparse} we will see that  the function $v_X$ might   as well  have  different  superior and  inferior exponential growth rates $\omega^{\pm} (X)$ (notice, in contrast, that $\delta(\Gamma)$  always is a true limit). 
 \vspace{3mm}

The main result of the paper concerns finiteness of the Bowen-Margulis measure and an aymptote for the volume growth function of $\frac14$-pinched spaces with lattices:
 
\begin{theorem}
\label{teor14}
Let $X$ be  a Hadamard space with   curvature   $-b^2\leq K_X \leq -a^2 $, and let $\Gamma$ be a nonuniform lattice of $X$. 
If   $\bar X = \Gamma \backslash X$  has {\it asymptotically $1/4$-pinched curvature} (that is,  for any $\epsilon >0$,   the metric satisfies 
 $-k_+^2\leq K_{X} \leq -k_-^2 $ with $k_+^2 \leq 4k_-^2 + \epsilon$  outside some compact set $\bar C_{\epsilon} \subset \bar X$),  then:

\noindent (i) $\Gamma$ is divergent and the  Bowen-Margulis measure $\mu_{BM}$   of  $U\bar X$   is finite;

\noindent (ii)   $\omega^{+}(X) = \omega^- (X) =\delta(\Gamma)$;

\noindent (iii)  there exists a function $\bar m(x) \in L^1 (\bar X)$  such that  $ v_X(x, R) \sim m(x) e^{\delta(\Gamma)R}$, where $m(x)$ is the lift of $\bar m$ to $X$.

\end{theorem}

\noindent From (i) it  follows  that {\em the geodesic flow of any asymptotically  $\frac14$-pinched, negatively curved  manifold  of finite volume  is ergodic  and totally conservative} w.r. to  $\mu_{BM}$, by   Hopf-Tsuji-Sullivan Theorem (see \cite{sullivan}, \cite{roblin}), contrary to the case  of general negatively curved  manifolds of finite volume (e.g., those obtained from convergent lattices).  \\
Condition (iii) also implies that {\em volume equidistributes on large spheres}, i.e. the volume $v_X^\Delta (x, R)$ of annuli in $X$  of thickness $\Delta$  satisfies the precise asymptotic law: 
\[ v_X^\Delta (x, R) \sim 2 m(x) \sinh(\Delta \delta(\Gamma)) e^{\omega(X)R} \]
Notice that the above theorem   also covers  the classical case of noncompact symmetric spaces of rank one (where the proof of the divergence and the asymptotics is direct). 
\pagebreak

One may wonder about the meaning   of the $\frac14$-pinching condition. This turns out to be an asymptotic, geometrical condition on the influence and wildness of maximal parabolic subgroups of $\Gamma$ associated to the cusps of $\bar X=\Gamma \backslash X$. Parabolic groups, being elementary, do not necessarily have a critical exponent which can be interpreted as  a true limit; rather, for a parabolic group of isometries $P$ of $X$, one can consider the limits

$$ \delta^+ (P) = \limsup_{R \rightarrow \infty} \frac{1}{R} \ln v_P (x,R) \, \hspace{1cm}
\delta^- (P) = \liminf_{R \rightarrow \infty} \frac{1}{R} \ln v_P (x,R) $$
and the critical exponent $\delta (P)$ of the Poincar\'e series of $P$ coincides with $\delta^+(P)$. 
\\ 
Accordingly, we say that a lattice $\Gamma$ is {\em sparse}  if it has a maximal parabolic subgroup $P$ such that 
$ \delta^+ (P)>  2  \delta^- (P) $ (conversely, we will say that $\Gamma$ is {\em parabolically $\frac12$-pinched} if it is not sparse).
Such parabolic groups in  $\Gamma$, together with dominant parabolic subgroups,  are precisely associated to cusps whose growth   can wildly change, and this can globally influence the growth function of $X$. Namely, we can prove:

\begin{theorem}
\label{teor14basic}
Let $X$ be  a Hadamard manifold with pinched, negative curvature  \linebreak $-b^2\leq K_X\leq -a^2 <0$.
If $X$ has a nonuniform lattice $\Gamma$ which is neither exotic nor sparse,  then $\Gamma$ is divergent  with finite Bowen-Margulis measure; moreover, $v_{X} \asymp v_{\Gamma}$ and  $X$ has a Margulis function $m(x)$, whose projection is $L^1$ on $\bar X = \Gamma \backslash X$. 
\end{theorem}

 The divergence and finiteness of the Bowen-Margulis measure in Theorem \ref{teor14} and Theorem \ref{teor14basic} are both   consequence of  a {\em critical gap}  between $\delta(\Gamma)$ and  the exponential growth rates   $\delta^+(P_i)$ of all parabolic subgroups; this will be proved in \S4. \linebreak 
In particular, we will see  that {\em any lattice $\Gamma$ in a  negatively curved, $\frac{1}{4}$-pinched space is never exotic} (nor sparse). 
For this, we will use an asymptotic characterization of  the hyperbolic  lattices  as the only lattices in  spaces $X$ with  pinched curvature   $-b^2 \leq K_X\leq -a^2$ realizing the least possible value for the topological entropy of $\bar X =\Gamma \backslash X$, i.e. satisfying $\delta(\Gamma)=(n-1)a$. 
In the compact case, this result can be deduced from Knieper's work on spherical means (following the proof of Theorem 5.2, \cite{knieper}), or from Bonk-Kleiner \cite{BK}  (for convex-cocompact groups);  on the other hand, see \cite{dpps-rigidity} for  a complete  proof in the case of non-uniform lattices and  the  analysis of  the new difficulties arising in the non-compact case. 

The existence of the Margulis function in Theorems \ref{teor14} and \ref{teor14basic}   relies  on a Counting Formula (Proposition \ref{propdissection}), proved in \S3; the formula enables us to reduce the computation of $v_X$ to the {\em analytic profile} of the cusps of $\bar X$ and $v_\Gamma$ (so, in the last instance,   to T.Roblin's asymptotics {\em (a)\&(b)}).
\vspace{3mm}

The last  part of the paper is devoted to studying sparse and exotic lattices, to understand the necessity of the $\frac14$-pinching (or $\frac12$-parabolically pinching) conditions. \linebreak
The following result shows that  Theorem \ref{teor14basic} is the best that we can expect for Hadamard spaces with quotients of finite volume:

\pagebreak 
\begin{theorem}
\label{teorex}
Let  $X$ be a Hadamard manifold with pinched negative  curvature  \linebreak $-b^2\leq K_X\leq -a^2 <0$ admitting  a nonuniform lattice $\Gamma$.
 \vspace{1mm}
        
\noindent {\bf (i)}  If $\Gamma$ is exotic and the dominant subgroups $P$ satisfy $ \delta  (\Gamma) = \delta^+ (P)  < 2  \delta^- (P) $, then both $v_X$ and $v_{\Gamma}$ are purely exponential or lower-exponential, with the same exponential growth rate $\omega (X)=\delta(\Gamma)$. Namely:
\vspace{-2mm}
        
\begin{itemize}[leftmargin=5mm]
\itemsep0mm
\item either $\mu_{BM}<\infty$,   and  then  $v_{X}$ is purely exponential and $X$ has a Margulis function;
   
\item or  $\mu_{BM} = \infty$, and in this case  $v_{X}$ is lower-esponential.
\end{itemize}

\vspace{-2mm}        
\noindent The two  cases can actually occur, cp. Examples \ref{exexotic<}(a)\&(b).
 
 \vspace{1mm}    
\noindent {\bf (ii)} If  $\Gamma$ is exotic and a dominant subgroup $P$ satisfies $ \delta  (\Gamma) =   \delta^+ (P) =2  \delta^- (P) $, then $\omega (X)=\delta(\Gamma)$ but in general $v_X \not \asymp v_{\Gamma}$, and $X$ does not admit a Margulis function. Namely, there exist cases  (Examples \ref{exexotic=}(a)\&(b)) where:
 \vspace{-2mm}           

\begin{itemize}[leftmargin=5mm]
\itemsep0mm
\item $\mu_{BM}<\infty$,   with  $v_{\Gamma}$ purely exponential and  $v_{X}$ upper-exponential;

\item  $\mu_{BM} = \infty$, with  $v_{\Gamma}$ lower-exponential and  $v_{X}$ upper-exponential.
\end{itemize}
\end{theorem}
By lower- (respectively, upper-) exponential, here, we mean a function $f$ with  exponential growth rate  $\omega=\limsup_{R \rightarrow \infty} \frac{1}{R} \ln f(R)$, but such that   
$\liminf_{R \rightarrow \infty} f(R)/e^{\omega R} =0$ (resp.  $\limsup_{R \rightarrow \infty}  f(R)/e^{\omega R} =+\infty$).
\vspace{2mm}        
 
 We shall see that all these examples can be obtained as lattices in $(\frac14 -\epsilon)$-pinched spaces, for arbitrary $\epsilon>0$, which shows the optimality of the 
$\frac14$-pinching condition.

\noindent On the other hand, if $\Gamma$ is sparse,  one can even have $\omega^+ (X)> \omega^- (X)> \delta(\Gamma)$, and the Example  \ref{exsparse} shows that virtually any asymptotic behaviour for $v_X$ can occur. 
Thus, the case of exotic lattices with a parabolic subgroup such that $\delta^+(P)=2\delta^-(P)$ can be seen as the critical threshold where a transition happens, from   functions $v_{\Gamma}, v_X$ with same asymptotic behaviour to functions with even different exponential growth rate.

\noindent Notice at last  that  the condition  $ \delta^+ (P)  < 2  \delta^- (P) $  is satisfied when  ${b^2\over  a^2}<{1\over 4}$, and that this last condition implies that  the group $P$ is abelian \cite{BK}.
\vspace{2mm}

\footnotesize 

\noindent {\sc Notations}.  Given two functions $f,g: \mathbb{R}_+ \rightarrow \mathbb{R}_+$, we  will systematically write  
 $f\stackrel{C}{\prec}g $ for $R>R_0$ (or $g \stackrel{C}{\succ} f$)  if   there exists $C>0$ such that  $f(R) \leq C g( R)$  for these values of $R$. We say  that $f$ and $g$ are {\em weakly asymptotically equivalent} and write $f  \stackrel{C}{\asymp} g$  when  $g  \stackrel{C}{\prec}  f \stackrel{C}{\prec}  g$ for $R \gg0$; we will simply write $f \asymp g$  and $f \prec g $   when the constants $C$ and $R_0$ are unessential. We say that $f$ and $g$ are {\em asymptotically equivalent} and write $f \sim g$ when  $\lim_{R\rightarrow +\infty}f (R) /g(R) =1$.
 
\noindent We define the upper and lower exponential growth rates of  a function $f $  respectively as: 
$$ \omega^+ (f) ={\limsup_{R\to+\infty}}  \;  R^{-1} \ln f(R)  
\quad \mbox{ and } \quad 
\omega^-(f)=\omega(f)={\liminf _{R\to+\infty}}  \;   R^{-1}\ln f(R) $$

\noindent and we simply write $\omega (f)$ when the two limits coincide. Also, we will say that   $f$ is {\em purely exponential} if $f \asymp e^{\omega(f)R}$, and that  $f$   is
{\em lower-exponential} (resp. {\em upper-exponential}) when $\liminf_{R\to+\infty} { f(R)) \over  e^{\omega(f)R}}=0$ 
(resp.  $\limsup_{R\to+\infty} {f (R) \over  e^{\omega(f)R} }=+\infty$).

\noindent Finally, if $f $ and $g $ are two real functions, we will use the notation 
$f\ast_{\!\Delta} g$ for the discrete convolution of $f$ and $g$ with gauge $\Delta$, defined by
$\displaystyle (f\ast_{\!\Delta} g)(R)=\sum_{h,k \geq 1}^{h+k=\lfloor R/\Delta \rfloor} f(h\Delta)g(k \Delta) $.
We notice here that, for nondecreasing functions $f$ and $g$, this is weakly equivalent to the usual convolution, namely
$$\Delta \cdot  (f\ast_{\Delta} g) \, (R-\Delta) \leq (f\ast  g) \, (R) = \int_0^R  f(t)g(R-t)dt 
   \leq  2 \Delta \cdot  ( f\ast_{\Delta} g) \, (R+2\Delta).
$$

 \pagebreak

 \normalsize

 \section{Growth of parabolic subgroups and of lattices modulo parabolic subgroups}
 
{\em \underline{Throughout all the paper, unless otherwise stated, $X$ will be a Hadamard space of} \linebreak \underline{dimension $n$, with pinched  negative sectional curvature $-b^2 \leq K_X \leq -a^2 <0$}.}

For  $x,y \in X$ and $\xi$ belonging to the geometric boundary $ X (\infty)$, we will denote $[x,y]$ (resp. $[x,\xi]$) the geodesic segment  from $x$ to $y$ (resp. the ray  from $x$ to $\xi$).\linebreak
We will  repeatedly make use of the following, classical result in strictly negative curvature: there exists $\epsilon (a, \vartheta)= \frac{1}{|a|}\log (\frac{2}{1 - \cos \vartheta})$ such that any geodesic triangle   $xyz$   in $X$ making  angle $\vartheta=\angle_z (x,y) $ at  $z$ satisfies: 
\begin{equation}
\label{eqoppositetriangle} d(x,y) \geq d(x,z) + d(z,x) - \epsilon(a, \vartheta).
\end{equation}



Let $b_{\xi}(x,y)=\lim_{z\to\xi}d(x,z)-d(z,y)$ be the Busemann function centered at $\xi$. \linebreak
The level set $\partial H_\xi (x)\!=\!\{ y \, | \, b_\xi (x,y)\!=\!0\}$ (resp. the suplevel set \nolinebreak $H_\xi (x)\!=\!\{ y \, | \, b_\xi (x,y) \!\geq\! 0\}$ is the horosphere (resp. the horoball) with center  $\xi$ and passing through $x$. \linebreak
 From (\ref{eqoppositetriangle})  we easily deduce the following:

\begin{lemma}
\label{lemmahor} 
 For any $d>0$, there exists $\epsilon_1=\epsilon_1(a,d) \geq \epsilon (a, \frac{\pi}{2})$  with the following property: 
given  two disjoint horoballs  $H_1, H_2$ at distance $d= d(H_1, H_2) =d (z_1,z_2)$ with $z_i \in \partial H_i$, then for any 
  $x \in H_1$ and $y \in H_2$ we have 
$$   d(x,z_1) + d(z_1, z_2) + d(z_2, y) - \epsilon_1 (a,d) \leq d(x,y) \leq   d(x,z_1) + d(z_1, z_2) + d(z_2, y).$$ 
\end{lemma}

\noindent  {\bf Proof.}   As  $K_X \leq -a^2$ and horoballs are convex, for any $y \in H_2$   the angle \linebreak
 $\vartheta (y) = \angle_{z_1} z_2, y$ satisfies $ \tan \vartheta (y) \leq \frac{1}{\sinh (d/|a|)}$
  (cp. for instance \cite{sam}, Prop.8).  
Then, we have   $\angle_{z_1} x, y \geq  \frac{\pi}{2} - \vartheta (y) \geq \vartheta(d)$ with $ \vartheta(d) >0$ for $d \neq 0$, hence,  by
(\ref{eqoppositetriangle}),
$$
d(x,y) \geq d(x,z_1) +  d(z_1y)  -  \epsilon (a ,   \vartheta(d)) \geq d(x,z_1) +  d(z_1,z_2) +  d(z_2,y)  - \epsilon_1(a,d) 
$$
for $ \epsilon_1(a,d) = \epsilon (a, \vartheta(d))+  \epsilon (a, \frac{\pi}{2})$.$\Box$

Let $d_\xi$  denote the horospherical distance between two points on a same horosphere centered at $\xi$. 
If $\psi_{\xi,t}: X \rightarrow X$ denotes the radial flow in the direction of $\xi$, we define:
\begin{equation}
\label{eqt}
t_\xi (x,y) = 
\left\{  \begin{array}{ll}
\inf \{t>0 \; | d_\xi (\psi_{\xi, t+ \Delta } (x) ,\psi_{\xi, t} (y)) <1\}  &  \mbox{ if } b_\xi (x,y)=\Delta \geq 0 ;   \\
\inf \{t >0 \; | d_\xi (\psi_{\xi, t} (x),\psi_{\xi,t-\Delta } (y))   <1\}  &    \mbox{ if } b_\xi (x,y)=\Delta < 0.     \end{array} \right.
\end{equation}
\noindent If $y$ is closer to $\xi$  than $x$, let $x_{\Delta}= [x,\xi[ \cap \partial H_{\xi}(y)$: then, $t_\xi (x,y)$ represents the minimal time we need to apply the radial flow $\psi_{\xi,t}$ to the points $x_{\Delta}$ and $y$ until they are at horospherical distance less than $1$. Using (\ref{eqoppositetriangle}) and the lower curvature bound $K_X \geq -b^2$,  we obtain in  \cite{DPPS} the following estimate, which is also crucial in our computations:

 \begin{applemma}
\label{lemmapp} 
${}$

\noindent There exists $\epsilon_0 =\epsilon_0 (a,b)\geq \epsilon (a, \frac{\pi}{2})$  such that   for all $x,y \in X$ and $\xi \in X (\infty)$ we have: 
\vspace{-3mm}        

$$  2t_\xi (x,y)  + |b_\xi (x,y)|  - \epsilon_0 \leq d(x,y) \leq   2t_\xi (x,y)  + |b_\xi (x,y)| + \epsilon_0$$ 
\end{applemma}


In this section we give estimates  for the growth of annuli in a parabolic subgroup  and in quotients of a lattice  by a parabolic subgroup, which will be used later. So,   let us fix some notations.  
We let  $A^\Delta(x,R)= B\left(x,R+\frac{\Delta}{2} \right) \setminus B \left(x,R-\frac{\Delta}{2} \right)$ be the annulus   of radius $R$ and thickness $\Delta$ around $x$. For a group $G$ of isometries of $X$, we will consider  the orbital functions
$$v_G (x,y,R) =  \# \left( B(x,R) \cap Gy \right)  \hspace{1cm}  
v^\Delta_G (x,y,R) =  \# \left( A^\Delta(x,R) \cap Gy \right)  $$
and we set $v_G (x,R) = v_G(x,x,R)$,  $v^\Delta_G (x,R) = v^\Delta_G(x,x,R)$ and $v^\Delta_G (x,R) = \emptyset$ for $\Delta<0$. 
We will also need to consider the growth function of   coset spaces, endowed with the natural quotient metric:  
 if $H<G$, we define $d_x (g_1 H, g_2 H):= d(g_1 Hx, g_2 Hx)$ and 
\vspace{-3mm}        

$$v_{G/H} (x,R) :=  \# \{ gH \; | \;   | gH |_x = d_x (H, gH)  < R \}$$

\vspace{-3mm}        

$$v^\Delta_{G/H} (x,R) = v_{G/H} \left(x,R+\frac{\Delta}{2}\right) -   v_{G/H} \left(x,R-\frac{\Delta}{2}\right) \,.$$ 

\noindent We will use analogous notations for the growth functions of  balls and annuli in  the spaces of left and double cosets $H\backslash G$, $ H \backslash  G \slash  H $ with the metrics

$$d_x (Hg_1,  Hg_2)  := d(H g_1 x,H g_2 x) = | g_1^{-1} H g_2 |_x $$ 

 \vspace{-3mm}       

$$d_x (Hg_1H,  Hg_2H)  := d(Hg_1H x ,Hg_2H x) =  | g_1^{-1} H g_2 H |_x \;.$$

 The growth of the orbital function of a bounded parabolic group $P$ is best expressed by introducing the horospherical area function.
 Let us recall the necessary  definitions:

\begin{definitions} 
\label{defiparabolic}
{\em 
\noindent Let $P$ be a {\em bounded} parabolic group of $X$ fixing $\xi \in X (\infty)$: that is,   acting cocompactly on $X (\infty) - \{ \xi \}$ (as well as on every horosphere $\partial H$ centered at $\xi$). \linebreak 
Given $x\in X$,  let ${\cal D} (P,x)$  be a  {\em Dirichlet  domain}  centered at $x$   for the action of $P$ on $X$; \linebreak that is, a convex  fundamental domain contained in the closed subset 
\vspace{-2mm}        

 $$ \overline {\cal D} (P,x) = \{ y \in X  \; | \; d(x,y) \leq d(px,y) \mbox{ for all } p \in P \}$$ 

\noindent We set 
  ${\cal S}_x ={\cal D} (P,x) \cap \partial H_\xi (x)$ and   ${\cal C}_x ={\cal D} (P,x) \cap H_\xi (x)$, and denote by  ${\cal S}_x (\infty)$  the trace at infinity of ${\cal D} (P,x)$, minus  $\xi$; 
  these are, respectively,   fundamental domains for the actions of  $P$ on  $\partial H_\xi (x), H (x)$ and $X(\infty) \!-\! \{\xi\}$. \\
The {\em horospherical area  function} of $P$ is the function
        \vspace{-3mm}

$${\cal A}_P (x,R) = \vol \left[ P  \backslash \psi_{\xi,R} \left(  {\partial H_{\xi} (x)}  \right) \right]  
= \vol \left[   \psi_{\xi,R} \left(  {\cal S}_x   \right) \right]  $$

\noindent where the $\vol$  is  the Riemannian measure of horospheres.   We also define the  {\em cuspidal function} of $P$, which is the function  
  \vspace{-3mm}      

$${\cal F}_P (x,R) = \vol \left[ B(x,R) \cap H_\xi (x) \right]$$
that is, the volume of the intersection of a ball centered at $x$ and the horoball centered at $\xi$ and passing through $x$.
Notice that the  functions  ${\cal A}_P (x,R), {\cal F}_P (x,R)  $ only depend on the choice of the initial horosphere $\partial H_\xi (x)$. 
}
\end{definitions}

 \begin{remark}
{\em 
Well-known estimates  of the differential of the radial flow (cp. \cite{HIH}) 
 yield, when $-b^2 \leq K_X \leq -a^2 <0$,
 \vspace{-3mm}       

\begin{equation}
\label{jacobi-norm}
e^{-bt}  \parallel \! v \!  \parallel \leq  \parallel  \! d\psi_{\xi, t}(v)  \!  \parallel \leq e^{-at} \parallel  \!  v  \! \parallel
\end{equation}
Therefore we deduce that, for any $\Delta > 0$,

\begin{equation}
\label{eqA(R+C)}
e^{-(n-1)b\Delta}\leq  \frac{ {\mathcal A}_P(x, R+\Delta) }{{\mathcal A}_P(x, R) } \leq e^{-(n-1)a\Delta}
\end{equation}
}
 \end{remark}

The following Propositions show   how the horospherical area  ${\cal A}_P$ and the cuspidal function ${\cal F}_P $  are   related to the orbital  function of $P$; they  refine and precise some estimates given in  \cite{DPPS} for $v_P (x,R)$.

\begin{proposition}
\label{propVp}
Let $P$ be a bounded parabolic group of $X$ fixing $\xi$, with  $diam ({\cal S}_x) \leq d$. \\
There exist  $C=C(n, a,b, d)$ and  $C'=C'(n, a,b, d;\Delta)$ such that:
\vspace{-5mm}        

\begin{equation}
\label{eqVp} 
v_P (x,y, R)    \;\; \stackrel{C}{\asymp}  \;\;  {\cal A}_P^{-1} \left(x, \frac{R+ b_\xi (x,y)}{2} \right) 
\;\;\;\; \forall R \geq   b_\xi(x,y) \!+R_0 \hspace{25mm}
\end{equation}
\begin{equation}
\label{eqVdeltap}  
v^{\Delta}_P (x,y, R)    \; \stackrel{C'}{\asymp}  \;  {\cal A}_P^{-1} \left(x, \frac{R+ b_\xi (x,y)}{2} \right) 
\;\;\;\;  \forall R \geq  b_\xi(x,y) \!+\! R_0 \; \mbox{ and } \;\forall \Delta > \Delta_0
\end{equation}

\noindent for explicit constants $R_0$ and  $\Delta_0$
only depending on $n,a,b,d$.
\end{proposition}

\begin{proposition}
\label{propFp}
Same assumptions as in Proposition \ref{propVp}. We have:
\vspace{-3mm}        

\begin{equation}
\label{eqFp}  
{\cal F}_P (x, R) 
 \;\; \stackrel{C}{\asymp} \;\; 
 \int_0^R \frac{   {\cal A}_P (x,t)    }{    {\cal A}_P \left(x,\frac{R+t}{2}\right)   } dt
 \;\;\;\; \forall R \geq R_0
 \end{equation}
\end{proposition}

 \begin{remark}
 \label{reminequalities}
{\em 

\noindent  More precisely, we will prove (and use later) that:

\noindent (i) $\;v_P (x,y, R) \stackrel{C}{\prec}  {\mathcal A}^{-1}_P \! \left( x,\! \frac{R+b_\xi (x,y)}{2} \right)$ 
 for all $R>0$;

\noindent (ii) $v^{\Delta}_P (x,y, R) \stackrel{C'}{\prec}   {\mathcal A}^{-1}_P \! \left( x, \!\frac{R+b_\xi (x,y)}{2} \right)$ for all $\Delta, R \! >\!0$;


  \noindent (iv)  $\;  {\cal F}_P (x, R) \;   \stackrel{C}{\prec}   \;
 \int_{_0}^{^R} \frac{   {\cal A}_P (x,t)    }{    {\cal A}_P \left(x,\frac{R+t}{2}\right)   } dt \;\;$  for  all  $R>0$.
}
\end{remark}

As a direct consequence of (\ref{eqFp})   and  (\ref{eqVp}) we have  (see also Corollary 3.5 in \cite{DPPS}):

\begin{corollary}
\label{cordeltaF}
Let $P$ be a bounded parabolic group of $X$. Then:
\vspace{-3mm}        

\begin{equation}
\label{eqdeltaF}
\delta^-(P) \leq \omega^- ({\cal F}_P) \leq \omega^+ ({\cal F}_P) \leq max \{ \delta^+(P), 2 (\delta^+(P) - \delta^-(P)) \}  
\end{equation}
\end{corollary}

\vspace{3mm}
  {\bf Proof of Proposition \ref{propVp}.}  
Since $v_P (x,y, R) = v_P (y,x, R)$ and ${\cal A}_P  (x, R) = {\cal A}_P   (y, R - b_\xi (x,y) )$, we can assume that $t =  b_{\xi}(x,y) \geq 0$.
 If $z \in \partial H_{\xi}(y)$ and $d(x,z)=R$, we know by Lemma \ref{lemmapp}  that  
 $ 2 t_\xi (x,z) + t -\epsilon_0 \leq d(x,z) \leq  2 t_\xi (x,z) + t +\epsilon_0$, so  $| t_\xi (x, z) -  \frac{R-t}{2} | \leq \epsilon_0/2$.
We deduce that 
$d_\xi \left(   \psi_{\xi, \frac{R+t+\epsilon_0}{2}} (x),\psi_{\xi, \frac{R-t+\epsilon_0}{2}} (z)  \right)   \leq 1$, so the set
 $ \psi_{\xi, \frac{R-t+\epsilon_0}{2}} \left(B (x, R)\cap \partial H_{\xi}(y) \right)$ is contained in the unitary ball $B^+$ of  the horosphere $\partial H_{\xi} ( x^+)$,  centered at    $x^+\!\!=\psi_{\xi, \frac{R+t+\epsilon_0}{2}} (x)$.
Similarly,   if $R \!>\! t +\epsilon_0$  then \nolinebreak $t_\xi (x, z) \!>\!0$, so 
$d_\xi \left(   \psi_{\xi, \frac{R +t-\epsilon_0}{2}} (x),\psi_{\xi, \frac{R-t-\epsilon_0}{2}} (z)  \right)\geq 1$, and the set
$\psi_{\xi, \frac{R-t-\epsilon_0}{2}} \left(B (x, R)\cap \partial H_{\xi}(y) \right)$ contains the unitary ball $B^-$ of  $\partial H_{\xi} ( x^-)$, 
centered at the point $x^-=\psi_{\xi, \frac{R+t-\epsilon_0}{2}} (x)$. \linebreak
We know that, by Gauss' equation,   the sectional curvature of   horospheres  of $X$ is between $ a^2-b^2 $ and $ 2b(b-a)$ (see, for instance, \cite{bk}, \S1.4); therefore, there exist positive constants $v^-=v^-( a, b)$ and $v^+=v^+( a, b)$  such that $\vol (B^+) < v^+$ and  $\vol (B^-) > v^-$. \linebreak
Now, let ${\cal S}_y = \psi_{\xi, t} ({\cal S}_x)$ be the fundamental domain for the action of $P$ on $\partial H_\xi (y)$  deduced from ${\cal S}_x$. There are  at least $v_{\mathcal P}(x, y, R-d)$ distinct fundamental domains $p{\cal S}_y$ included in $B (x, R)\cap \partial H_\xi (y)$;  since the radial  flow $\psi_{\xi,  t}$ is equivariant with respect to the action of $P$ on the horospheres centered at $\xi$, there are also at least $v_P(x, y, R-d)$ distinct fundamental domains $\psi_{\xi, \frac{R-t+\epsilon_0}{2}} (p{\cal S}_y)$ included in 
$\psi_{\xi, \frac{R-t+\epsilon_0}{2}}(B (x, R)\cap \partial H_{\xi}(y))$.  \linebreak
We deduce that
$ v_P (x, y,R-d) \cdot {\mathcal A}_P (x,  \frac{R+t+\epsilon_0}{2}) <  v^+$   and, by  (\ref{eqA(R+C)}), this gives  
$ v_P (x, y, R ) \stackrel{C}{\prec}  {\mathcal A}^{-1}_P (x, {R  +t\over 2})$ for all $R \geq 0$.
On the other hand, if $R > t + \epsilon_0$, we can cover the set $B (x, R)\cap\partial  H_\xi(y)$ with 
$v_P (x, y, R+d)$  fundamental domains $p{\cal S}_y$, with $p\in P$;  then,  again,
$\psi_{\xi, \frac{R-t-\epsilon_0}{2}}(B (x, R)\cap \partial H_{\xi}(y))$ can be covered by $ v_P(x, y,R+d)$ fundamental  
domains 
$\psi_{\xi, \frac{R-t-\epsilon_0}{2}} (p{\cal S}_y)$  as well, hence we deduce that
$v_P (x, y,R+d ) \cdot  {\mathcal A}_P (x, {R+t-\epsilon_0\over 2}) \geq  v^- $. This  implies that
$v_P (x, y,R) \stackrel{C}{\succ}    {\mathcal A}_P^{-1} (x, {R + t \over 2})$ for all $R>t+R_0$,  for $R_0=\epsilon_0 +d$ and a constant $C=C(n,a,b,d)$.

\noindent To prove the weak equivalence (\ref{eqVdeltap}), we just write, for $R+\frac{\Delta}{2}> t+R_0$:
$$     v^{\Delta}_P (x,y,R) = v^{\Delta}_P (x,y,R+ \Delta/2 ) - v^{\Delta}_P (x,y,R- \Delta/2 )
\geq \frac{  C^{-1}  }{  {\cal A}_P\left( \frac{R+t+\Delta/2}{2} \right)} - \frac{  C  }{  {\cal A}_P\left( \frac{R+t-\Delta/2}{2} \right)} $$

\vspace{-6mm}        

$$ \hspace{20mm} \geq   \frac{ C^{-1}    e^{(n-1)a\frac{\Delta}{4}} - C   e^{-(n-1)a\frac{\Delta}{4}}   }{    {\mathcal A}_P \left( x,\frac{R+t}{2} \right) } 
=   2 \sinh \left[ \frac14 (n-1)a\Delta  - \ln C\right] \cdot  {\cal A}^{-1}_P\left( \frac{R+ t}{2} \right)   $$

\noindent again by  (\ref{eqA(R+C)}),  if $\Delta >  \Delta_0=\frac{4\ln C}{(n-1)a}$.
Reciprocally, we have for all $R, \Delta> 0$:
$$       v^{\Delta}_P(x,y,R)  \leq  v_P(x,y,R+\frac{\Delta}{2})  
\leq \frac{  C  }{    {\mathcal A}_P \left( x,\frac{R+t+\Delta/2}{2} \right)   } 
\leq \frac{  C'(n,a,b,d;\Delta) }{  {\mathcal A}_P \left( x,\frac{R+t}{2} \right)    } \;\; \Box$$

 \noindent  {\bf Proof of Proposition \ref{propFp}.}  We  just integrate (\ref{eqVp})  over a fundamental domain ${\cal C}_x$ for the action of $P$ on $H_\xi (x)$:

$\displaystyle {\cal F}_P (x, R)   
   \! = \!  \sum_{p \in P} \vol [ B(x,R) \cap p{\cal C}_x] 
   \! = \!   \int_{{\cal C}_x} \sum_{p \in P}   \mathrm{1}_{B(x,R)} (pz) \; dz
   \! = \!   \int_{{\cal C}_x} v_P (x,y,R) \; dy $

%
  
\noindent so, integrating over each slice $\psi_{\xi, t}({\cal S}_x)$ by the coarea formula,   we obtain

$\displaystyle  
   \int_0^{R-R_0}  \!\!  \int_{ \psi_{\xi, t} ({\cal S}_x) }  \!\!   {\cal A}^{-1}_P \left(x, \frac{R+t}{2} \right)   dt 
 \; \stackrel{C}{\prec} \;  {\cal F}_P (x, R)  \;
  \stackrel{C}{\prec}  \; \int_0^R  \!\! \int_{ \psi_{\xi, t} ({\cal S}_x) }  \!\!  {\cal A}^{-1}_P \left(x, \frac{R+t}{2} \right)  dt$

\noindent (the left inequality holding for $R>R_0$). By (\ref{eqA(R+C)}),  both sides are weakly equivalent  to  the integral $\displaystyle \int_0^R \frac{ {\cal A}_P (x,t) }{ {\cal A}_P (x, \frac{R+t}{2}) } dt$,  up to a multiplicative  constant $c=c(n, a, b,d)$.$\Box$

%
%
%
%
%

\begin{remark}{\em 
Thus, we see that  the  curvature bounds imply that $ v^{\Delta}_P(x, R) \asymp v_P  (x, R) $ for $\Delta$ and $R$ large enough.
This also holds in  general for {\em non-elementary}  groups $\Gamma$ with finite Bowen-Margulis measure, as in this case
$v_\Gamma^{\Delta} (x,R) \sim  \frac{2 \parallel \! \mu_x \! \parallel^2}{\parallel \! \mu_{BM} \! \parallel}   \sinh [ \frac{\Delta}{2} \delta(\Gamma)] e^{\delta(\Gamma)R} $
by  Roblin's asymptotics.
On the other hand,  it is unclear whether  the weak  equivalence $ v^{\Delta}_\Gamma  \asymp v_\Gamma  $ holds for non-elementary  lattices $\Gamma$, when   $\parallel \!  \mu_{BM}  \!\parallel =\infty$. 
}
\end{remark}

In the next section we will also need  estimates for  the growth of annuli in the spaces of left and right cosets of a lattice $\Gamma$ of $X$, modulo a bounded parabolic subgroup  $P$.  \linebreak
 Notice that, if $P$ fixes $\xi \in X (\infty)$,   the function $v_{P  \backslash \Gamma} (x,R)$ counts the number of points   $\gamma x \in \Gamma x$  falling in the Dirichlet domain ${\cal D} (P,x)$ of $P$   with $d(x, \gamma x)<R$; on the other hand, the  function $v_{ \Gamma/ P } (x,R)$ counts the number of  horoballs  $\gamma H_\xi (x)$ at distance \linebreak (almost) less than $R$ from $x$. It is remarkable that, even  if these functions count geometrically distinct objects,  they are weakly asymptotically equivalent, as the following Proposition will show.  
 Actually, let   $H_\xi$ be  a horoball  centered at the parabolic fixed point  $\xi$ of $P < \Gamma$;  we call $depth(H_\xi)$   the minimal distance $\min_{\Gamma - P}  d(H_\xi, \gamma H_\xi)$. \linebreak
  Then, for ${\cal S} _x$ defined as in  Definition  \ref{defiparabolic} we have:

\begin{proposition}
\label{propGamma/P}
Let $\Gamma$ be a torsionless, non-elementary,  discrete group of isometries of $X$,  let $P$  a  bounded parabolic subgroup of $\Gamma$, and  let $x \in X$  be fixed.  
Assume that  $\max\{diam ({\cal S} _x), 1/ depth (H_\xi (x)) \} \leq d$,  
and let $\ell$ be the minimal displacement  $d(x, \gamma x) $ of the elements $\gamma \in \Gamma$ whose  domains of attraction   
${\cal U}^{\pm} (\gamma, x ) = \{ y \; | \; d (\gamma^{\pm 1} x, y) \leq  d (x, y) \}$ 
are  included in  the Dirichlet domain  ${\cal D} (P,x)$.

\noindent There exists a  constant  $\delta_0=\delta_0(a, d)$ such that,   for all  $\Delta, R>0$:

\noindent (i)  
$\;\;\; v^{\Delta-\delta_0}_{P \backslash \Gamma}  (x,R)  
\; \leq \; v^\Delta_{\Gamma \slash P} (x,R) 
\; \leq \;    v^{\Delta+\delta_0}_{P \backslash \Gamma} (x,R) $; 

\noindent (ii) 
 $\;\;\,   \frac12  v^{\Delta-2\ell}_{ \Gamma} (x,R)   
 \; \leq \; v^\Delta_{P \backslash \Gamma}  (x,R)  
 \; \leq \; v^\Delta_{\Gamma} (x,R)  $; 
 
\noindent (iii)  
$\;\,  
 \frac12   v^{\Delta-\delta_0-2\ell}_{\Gamma}  (x,R)  
\; \leq \; v^\Delta_{   \Gamma  \slash P  } (x,R) 
\; \leq \;    v^{\Delta +\delta_0}_{\Gamma} (x,R) $; 
 
\noindent (iv)  
$\;\;  
 \frac14  v_{\Gamma}^{\Delta -\delta_0 - 4\ell} (x,R)  
\; \leq \; v_{  P  \backslash  \Gamma  \slash P  }^{\Delta} (x,R) 
\; \leq \;    v_{\Gamma}^{\Delta} (x,R) $.
\end{proposition}

\noindent 
 
Notice that (iv) strenghtens a result  of S. Hersonsky and F. Paulin on the number of rational lines with depth smaller than $R$  (cp.  \cite{HP} Theorem 1.2,   where the authors furthermore assume the condition $\delta_P < \delta_\Gamma$). Actually, let $H_\xi$ be the largest  horosphere centered at $\xi$ non intersecting any other $\gamma    H_\xi $ for $\gamma \neq e$, and recall that the {\em depth of  a geodesic}  $c=(\xi, \gamma \xi)$ is defined as the length of the maximal subsegment $\hat c \subset c$ outside $\Gamma  H_\xi$. 
The double coset space $P  \backslash \left(\Gamma \!\! -\!\! P  \right) \slash P $ can be identified with the set of oriented geodesics $(\xi, \gamma \xi)$ of $X$ with $\gamma \in \Gamma  \!\!-\!\! P$.
Then, if $x \in \partial H_\xi $,  the counting function $v^\Delta_{P \backslash \left(\Gamma -  P  \right) \slash P } (x,R)$  corresponds to the number of geodesics of $\bar X = \Gamma \backslash X$ which travel a time about $R$ outside the cusp  $\bar {\cal C}  = P \backslash   H_\xi$, before entering and definitely staying (in the future and in the past) in  $\bar {\cal C}$. 
\vspace{3mm}

 {\bf Proof.}  The right-hand inequalities in (ii), (iii), (iv) are trivial. \\ Let us prove (i).  
  We first define two sections of the projections $ P\backslash \Gamma \leftarrow  \Gamma \rightarrow  \Gamma \slash P$. \linebreak
 Consider the fundamental domain ${\cal S}_x (\infty) $ for the action of $P$ on  $X(\infty) \!-\! \{\xi\}$ given in  \ref{defiparabolic}, and choose  for each $\gamma \in \Gamma$, a representative $\hat \gamma$ of $\gamma P$ which minimizes the distance to \nolinebreak $x$. \linebreak

 \noindent Then, we set 
  $$ \widehat  \Gamma = \{ \widehat  \gamma \; | \; \gamma P  \in \Gamma /P \}$$
 $$  \Gamma_0 = \{  \gamma_0 \; | \;    \gamma_0 \in \Gamma,   \gamma_0 \xi \in {\cal S}_x (\infty)  \} \cup \{ e\}.$$
 
\noindent We have bijections  $\widehat  \Gamma \cong  \Gamma /P$ and $\Gamma_0 \cong P \backslash \Gamma$, 
 as ${\cal S}_x (\infty)$ is a fundamental domain.
Moreover,  every $ \gamma_0 \in   \Gamma_0$  {\em almost} minimizes the distance to $x$ in its right coset $P\gamma_0$.
Actually,    for all $\gamma \in \Gamma$ 
   set
  $z({\gamma} ) = (\xi,  \gamma \xi)  \cap \partial H_{\xi} (x)$ 
  and  $z'( {\gamma}) = (\xi,   \gamma \xi)  \cap \gamma \partial H_{\xi} (x)$; then, 
  for all $p \in P$ we have, by Lemma \ref{lemmahor} 
\vspace{-3mm}

\begin{equation}
\label{eqalmostmin}
d( x, p \gamma_0  x) 
\geq d(x, p z({\gamma} )  ) +   d (pz({\gamma} ), pz'({\gamma} )) + d(pz'({\gamma} ), p   \gamma_0 x)  - \epsilon_1 (a,d)
\geq d(x,    \gamma_0 x) -c 
\end{equation}

 \noindent  as   $d(H_\xi (x), p   \gamma_0 H_\xi (x)) = d (pz({\gamma} ),pz'({\gamma} ))$, for $c   = 2d +  \epsilon_1 (a,d)$. \\
 We will now define a bijection  between  pointed metric spaces $i:   (P \backslash \Gamma, x_0)  \rightarrow   (\Gamma / P, x_0) $ 
which almost-preserves the distance to their base point $x_0=P$ 
(with respect to their quotient distances $ |\, \cdot \,  |_x = d_x (P, \cdot )$   as seen at the beginning of the section), as follows. \\
For every $ \gamma \in \Gamma$ we can write  $\gamma = \widehat \gamma p_\gamma$,  for uniquely determined $\widehat  \gamma \in \widehat  \Gamma$ and $p_\gamma \in P$; given a right coset $P \gamma$, we take $\gamma_0 \in \Gamma_0$ representing $P \gamma$ and then  set $i(P \gamma  ) := p_{\gamma_0 } \widehat   \gamma_0    P $.\linebreak 
The map $i$ is surjective. Actually, given $\gamma P$, we take $p \in P$ such that $p \gamma \xi  \in  {\cal S}_x (\infty)$, so that  $P \gamma   = P \gamma_0$, for $\gamma_0 = p \gamma \in \Gamma_0$; 
then, we write   $\gamma_0=\widehat{\gamma}_0 p_{\gamma_0}$, and  we deduce that \linebreak
 $i( P \gamma  ) =i( P \gamma_0  ) =  i( P \widehat{\gamma}_0 p^{-1} )   = p^{-1}   \widehat{\gamma}_0 P =  p^{-1}  \gamma_0  p_{\gamma_0}^{-1} P  = \gamma P$.\\
We now check that $i$ is injective. 
Given $\gamma_0 = \widehat{\gamma}_0 p_{\gamma_0}$ and $\gamma'_0 = \widehat{\gamma'}_0  p_{\gamma'_0 }  $ in $ \Gamma_0$ representing two right cosets $P \gamma$ and  $P \gamma'$, assume that 
 $p_{\gamma_0 } \widehat   \gamma_0    P = p_{\gamma'_0 } \widehat{\gamma'}_0    P $.
 Then,    $\widehat   \gamma_0 \xi = p   \widehat{\gamma'}_0 \xi$  for   $p = p_{\gamma_0 }^{-1}  p_{\gamma'_0 } \in P$, which yields $ p_{\gamma_0 } =  p_{\gamma'_0 }$  as   $\widehat   \gamma_0 \xi, \widehat{\gamma'}_0 \xi \in {\cal S}_x (\infty)$  and  $ {\cal S}_x (\infty)$ is a fundamental domain for the left action of $P$; so,  $\widehat   \gamma_0 P =   \widehat{\gamma'_0} P$, which implies that  $\widehat   \gamma_0 = \widehat{\gamma'_0} $ too (as $\widehat \Gamma$ is a section of $\Gamma / P$).
 Therefore,  $P\gamma = P \gamma_0 = P  \widehat{\gamma}_0 p_{\gamma_0} = P  \widehat{\gamma'}_0 p_{\gamma'_0} =  P \gamma'_0 = P \gamma'$.    \\
 To show that   $i$  almost preserves $| \; |_x $, we notice that, given a class $P\gamma$ and writing  its representative in $\Gamma_0 $ as $\gamma_0  = \widehat \gamma_0 p_{\gamma_0} $, we have
 \vspace{-3mm}       

 $$  | P \gamma |_x \leq  | \gamma_0  |_x \leq d(x, \widehat \gamma_0  x ) +  d( \widehat \gamma_0  x , \widehat \gamma_0  p_{\gamma_0 } x) 
 = | \widehat \gamma_0    |_x + | p_{\gamma_0 } |_x$$

\noindent while, by (\ref{eqalmostmin}) and by Lemma \ref{lemmahor}  
 $$| P \gamma |_x \geq | \gamma_0 |_x  - c
 \geq d(x, z'(\gamma_0  )) +  d(  z'(\gamma_0 ), \widehat \gamma_0    p_{\gamma_0 } ) - \epsilon_1(a,d) -c
  \geq | \widehat \gamma_0  |_x +  | p_{ \gamma_0 } |_x - 2c 
  $$
  
\noindent  as $d( z'(\gamma_0  ), \widehat \gamma_0 x) <d$.
On the other hand
$$ | i (P \gamma) |_x   
= | p_{\gamma_0} \widehat  \gamma_0  P |_x
\leq d (x, p_{\gamma_0} x) + d (   p_{\gamma_0}  x,    p_{\gamma_0}  \widehat \gamma_0  P  x) 
= |  p_{\gamma_0}  |_x + |  \widehat \gamma_0  |_x $$

\noindent while, as 
$z( p_{\gamma_0}  \widehat \gamma_0 )=   p_{\gamma_0}  z( \widehat \gamma_0 ) $ 
  and 
$z' ( p_{\gamma_0} \widehat \gamma_0 )=   p_{\gamma_0}  z'( \widehat \gamma_0 ) $, we get by Lemma \ref{lemmahor} 
$$| i (P\gamma) |_x 
\geq d (x, p_{\gamma_0} z ( \widehat  \gamma_0)  ) +  d ( p_{\gamma_0} z ( \widehat  \gamma_0) , p_{\gamma_0} \widehat  \gamma_0 P x) - \epsilon_1(a,d)   
 \geq | p_{\gamma_0} |_x + | \widehat \gamma_0 |_x - c.$$

\noindent This shows that  $ |   P\gamma   |_x  -c \leq | i( P \gamma) |_x \leq |   P\gamma   |_x   +2c$.
 \noindent We then immediately deduce 
 that $v_{P \backslash \Gamma}  (x,R-2c)  
\; \leq \; v_{\Gamma \slash P} (x,R) 
\; \leq \;    v_{P \backslash \Gamma} (x,R+c) $, as well as (i)  for $\delta_0= 4c $.
 
\pagebreak        
 \noindent The proof of the left-hand  inequality    in (ii) is a variation for annuli of a trick  due to  Roblin, cp. \cite{roblin}. 
 Actually, as $L(P) \subsetneq L(\Gamma)$, we can choose a   $\bar \gamma  \in \Gamma$ with $d(x, \bar \gamma x) =\ell $ and such that the domains of attraction 
${\cal U}^{\pm} (\bar \gamma, x )$  are   included in the   domain ${\cal D}(P,x)$.  \linebreak
\noindent Let $v_{{\cal D}(P,x)} (x,R)$ be the number of points of the orbit $\Gamma x$ falling in ${\cal D}(P,x) \cap B(x,R)$.         

\noindent We have:
$$v^\Delta_\Gamma (x, R) \leq v^\Delta_{{\cal D}(P,x)} (x, R) + v^{\Delta+2\ell }_{{\cal D}(P,x)} (x, R) \leq 2 v^{\Delta+2\ell }_{{\cal D}(P,x)} (x, R)   $$
since, for $\gamma x \in A^\Delta (x,R)$,   either $\gamma x \in  {\cal D}(P,x)$, or  $\bar \gamma  \gamma x \in {\cal D}(P,x)$ and 
$\bar \gamma  \gamma x \in A^{\Delta+2\ell} (x,R)$.
As the points of $P$ falling in ${\cal D}(P,x)$ minimize the distance to $x$ modulo the left action of $P$, we also have 
 $v^{\Delta+2\ell}_{{\cal D}(P,x)} (x, R) = v^{\Delta+2\ell}_{P \backslash \Gamma } (x,R)$, which  proves (ii).\\
 Assertion (iii) follows directly from (i) and (ii). 
  To show (iv), we need to estimate the number of classes $\gamma P$ modulo the left action of $P$, that is  the elements of $\widehat \Gamma$  such that $\widehat \gamma  x$ belongs to   the fundamental domain ${\cal D}(P,x)$. We choose an element  $\bar \gamma  \in \Gamma$ with ${\cal U}^{\pm} (\bar \gamma, x ) \subset {\cal D}(P,x)$ as before, and apply again Roblin's trick  to the classes $\gamma P$. 
 The set $\widehat \Gamma x$ can be parted in two disjoint subsets: the subset 
 $\widehat \Gamma_1 := \widehat \Gamma \cap   {\cal D}(P,x)$, and  the subset 
 $\widehat \Gamma_2 :=  \widehat \Gamma  \cap   {\cal D}(P,x)^c$,  whose elements $\widehat \gamma$ then satisfy
 $\bar \gamma \widehat \gamma \in  {\cal D}(P,x)$ and 
 $| \bar \gamma \widehat \gamma |_x \leq |  \widehat \gamma |_x + \ell$.\linebreak
 Then 
 $v^{\Delta}_{\Gamma / P}  (x,R)  
 = v^{\Delta}_{\widehat \Gamma_1}  (x,R) + v^{\Delta}_{\widehat \Gamma_2}  (x,R)
 \leq  2 v^{\Delta+2\ell}_{  P  \backslash  \Gamma  \slash P  } (x,R)$
 and we conclude by (iii).$\Box$


\section{Orbit-counting estimates for lattices}
 
In this section we give estimates
 of the orbital function  $v_\Gamma(x,y,R)$ and of $v_X (R)$  in terms of the orbital function of the parabolic subgroups $P_i$ and the associated cuspidal functions ${\cal F}_{P_i}$ of  $\Gamma$. These estimates will be used in \S\ref{margulis} and \S\ref{examples}; they stem from an accurate dissection of large balls  in compact and horospherical parts, assuming that  ambient  space $X$ admits a nonuniform lattice action.
        
\vspace{3mm}
Let $\Gamma$ be a lattice of $X$. The quotient manifold $\bar X= \Gamma \backslash  X $ is geometrically finite, and we have the following classical results  due to B. Bowditch \cite{bow}  concerning the structure of the  limit set $L(\Gamma)$ and of  $\bar X$:

  (a)  $L(\Gamma)= X (\infty)$ and it is the disjoint 
union of the radial limit set $L_{rad}(\Gamma)$ with  finitely many orbits
$ L_{bp} \Gamma =  \Gamma \xi_{1} \cup  \ldots \cup \Gamma \xi_{m}$ of {\em bounded} parabolic fixed 
points; this means that each  $\xi_i \in L_{bp}G $  is the fixed point of some  maximal bounded parabolic subgroup  $P_i$ of $\Gamma$;

  (b) {\em (Margulis' lemma)}  there exist  closed horoballs $H_{\xi_{1}}, \ldots, H_{\xi_{m}}$ centered respectively  at $\xi_{1}, \ldots, \xi_{m}$, such that  $gH_{\xi_{i}} \cap H_{\xi_{j}} = \emptyset$  for  all $1\leq i,j\leq m$ and  all $\gamma \in \Gamma - P_i$;

 (c)   $\bar X$  can  be decomposed into
 a disjoint union of a compact set $\bar  {\cal K}$ and finitely many ``cusps'' $\bar {\cal C}_1, ..., \bar {\cal C}_m$: each $\bar {\cal C}_i$ is isometric to the quotient of $H_{\xi_{i}}$  by  the maximal bounded parabolic group $P_i \subset \Gamma$. We refer to $\bar  {\cal K}$
and to  $\bar {\cal C} = \cup_i \bar {\cal C}_i$ as to the {\em compact core}  and the {\em cuspidal part} of $\bar X$.

Throughout this section, we  fix  $x \in X$
and  we consider a Dirichlet  domain ${\cal D}(\Gamma,x)$  centered at $x$; this is a   convex fundamental subset, and we may assume that  ${\cal D}$ contains the geodesic rays $[x,\xi_i[$.  Accordingly, setting 
 ${\cal S}_i = {\cal D} \cap \partial H_{\xi_i}$ and
 ${\cal C}_i = {\cal D} \cap H_{\xi_i} \simeq   {\cal S}_i \times \mathbb{R}_+$,
the fundamental  domain ${\cal D}$ 
 can be decomposed into a disjoint union:

$${\cal D} =  {\cal K}  \cup {\cal C}_1 \cup \cdots \cup {\cal C}_m$$

\noindent where  $ {\cal K}$ is a convex, relatively compact set containing $x$ in its interior  (projecting to a subset $\bar  {\cal K}$ in $\bar X$), while   ${\cal C}_i$  and  ${\cal S}_i$  are, respectively, connected fundamental domains for the action of $P_i$ on $H_{\xi_{i}}$ and  $\partial H_{\xi_{i}}$ (projecting respectively to subsets  $\bar {\cal C}_i$,   $\bar {\cal S}_i$ of  $\bar X$).

\noindent  Finally, as $L(P_i) =\{ \xi_i\}$, for every $1\leq i \leq m$ we can find an element $\gamma_i \in \Gamma$,  with $\ell_i= d(x, \gamma_i x ) $,  which is {\em in Schottky position with $P_i$ relatively to $x$}, i.e. such that the domains of attraction 
 ${\cal U}^{\pm} (\gamma_i) = \{ y \; | \; d (\gamma_i^{\pm 1} x, y) \leq d (x, y) \}$ 
are  included in the Dirichlet domain ${\cal D}(P_i,x)$, as in Proposition \ref{propGamma/P}.

\noindent {\em For the following,  we will then set
%
$d   = max \{ diam(  {\cal K}), diam (  {\cal S}_i), 1/depth (H_{\xi_i}), \ell_i \} \geq \epsilon_0.$}

%


\begin{proposition}[Counting Formula]
\label{propdissection}
${}$ \\
Assume that  $x,y \in X$ project respectively to the compact core $\bar  {\cal K}$ and to a  cuspidal end $\bar {\cal C}_i$ of $\bar X = \Gamma \backslash X$. There exists  $C''=C(n, a,b, d)$  such that:  
\vspace{-5mm}

%
 
$$
  \left[ v_\Gamma (x, \cdot)  \ast  \!  v_{P_i} (x,y,\cdot) \right]  \! ( R \! -\! D_0)
   \stackrel{C''}{\prec}   v_\Gamma (x,y,R)   \stackrel{C''}{\prec}  
  \left[ v_\Gamma (x, \cdot )  \ast   \!  v_{P_i} (x,y,\cdot)\right] \! (R \! +\! D_0)
   \;\;\; \forall R \! \geq \! 0
$$

%
\noindent for a constant $D_0$ only depending  on $n,a,b,d$. 
\end{proposition}

\noindent  {\bf Proof.}   We will write, as usual,   $|\gamma|_x=d(x,\gamma x)$ and    $|\gamma P|_x=d(x,\gamma Px)$, and   choose a constant  $\Delta > \max\{R_0,\Delta_0, 2 \delta_0 +4d\}$, where  $R_0, \Delta_0, \delta_0$ are the constants of Propositions \ref{propVp} and \ref{propGamma/P}. We first show that 
  \vspace{-3mm}       
 
\begin{equation}
\label{eqdissection1}
B(x, R) \cap \Gamma y   \;\; \subset 
 \hspace{-1mm}   \bigcup^{ N }_{ \scriptsize \begin{array}{c} k  = 1 \end{array}}
 \hspace{-1mm}   \bigcup_{ \scriptsize \begin{array}{c} \bar \gamma \in  \Gamma,     |\bar \gamma|  \leq k\Delta     \end{array}}
 \hspace{-5mm}  
 B \left(\bar \gamma x, (N-k) \Delta  \right)   \cap   \left( \bar \gamma P_i \right) \! y
\end{equation}

\noindent for $N=\lfloor \frac{R}{\Delta} \rfloor +2$.  
 Actually, let $\gamma y \in B(x, R) \cap  \gamma H_{\xi_i}$  and set $\bar y_i = [x,  \gamma \xi] \, \cap\,  \gamma \partial H_{\xi_i}$. \linebreak
By using the action of the group $\gamma P_i \gamma^{-1}$ on $\gamma H_{\xi_i}$, we can find  $\bar \gamma =  \gamma p$, with $p \in P_i$,  such that   $\bar y_i \in \bar \gamma {\cal  C}_i$.
Since  the angle  $\angle_{\bar y_i} (x,   \gamma y)$ at   $\bar y_i$ is greater than $\frac{\pi}{2}$, we have: 
\vspace{-3mm}        

$$ d (x, \gamma y) \leq  d( x, \bar y_i) +  d( \bar y_i,  \gamma y)
      \leq d( x,  \gamma y) + \epsilon_0  <  R+  \epsilon_0$$
   
\noindent with $ | \bar   \gamma |  \leq   d(x,\bar y_i) +d < R+d + \epsilon_0 \leq N \Delta$.
Then, if $k\Delta \leq | \bar \gamma | < (k+1)   \Delta$, we deduce

 $$ d(\bar \gamma x, \gamma y ) 
   \leq d( \bar y_i , \gamma y ) + d \leq R + \epsilon_0 - d(x, \bar \gamma x) +2d < (N-k) \Delta$$

\noindent which shows that    
$\gamma y = \bar \gamma p^{-1} y \in B (\bar \gamma x, (N-k) \Delta  ) \cap  (\bar \gamma P_i) y 
  = \bar \gamma \left[ B (x, (N-k) \Delta  ) \cap  P_i y \right]$.  \\
 Thus,  we obtain:
 \vspace{-5mm}
 
$$v_\Gamma (x,y,R) \leq
    \sum_{k= 1}^{N}   v_{\Gamma} \! \left(x, k \Delta \right) \cdot v_{P_i} \! \left(x,y, (N-k) \Delta \right) 
    \prec v_\Gamma \ast v_{P_i} (R + 2\Delta)$$
 

\noindent This proves the right hand side  of our inequality. \\
The left hand is more delicate, as we need to dissect the ball $B(x,R)$ in disjoint annuli.   So, consider the set  $\widehat \Gamma_i$ of minimal representatives  of $\Gamma / P_i$ as in the proof of Proposition 
\ref{propGamma/P}. We have: 

\begin{equation}
\label{eqdissection2}
A^{4\Delta}(x, R ) \cap \Gamma y  \;\; \supset
  \hspace{2mm} \bigsqcup_{k =0}^{N}    
    \hspace{-2mm}   \bigsqcup_{ \scriptsize \begin{array}{c}   \widehat \gamma  \in  \widehat \Gamma_i  
                                              \\     k\Delta  -\frac{\Delta}{2} \leq |  \widehat  \gamma  |  <  k\Delta  +\frac{\Delta}{2}   \end{array}}
  \hspace{-15mm}     A^{ \Delta } \left(\widehat \gamma x,   (N-k) \Delta  \right)   \cap   \left( \widehat \gamma P_i \right) \!y
\end{equation}

\noindent for $N= \lfloor \frac{R}{\Delta} \rfloor +1$. In fact, given  $  \gamma  y = \widehat \gamma p_i y \in A^{\Delta} \!\left(\widehat \gamma x, (N \!-\!k)\Delta \right)$ with 
$\widehat \gamma x \in A^{\Delta} (x, k\Delta )$ 
we  have again
\vspace{2mm}        

\noindent \hspace{4mm}  
$       N  \Delta - 2\Delta  
\leq    | \widehat \gamma|  + d (\widehat \gamma x, \gamma y) -2d - \epsilon_0 
\leq    d(x,   \gamma y ) 
\leq    |\widehat \gamma| + d (\widehat \gamma x, \gamma y) 
<  N \Delta +\Delta$

 \vspace{1mm}       
\noindent as $\Delta > 2d + \epsilon_0$, hence $\gamma y \in A^{4\Delta} (x,R )$.
Notice that (\ref{eqdissection2}) is a disjoint union, as the annuli with the same center do not intersect by definition, while for $\widehat \gamma \neq \widehat \gamma'$ the orbits   $\widehat \gamma P_i y$ and  $ \widehat \gamma'  P_i y$ lie  on different horospheres $\widehat \gamma H_i \neq \widehat \gamma' H_i$, which are disjoint  by Margulis' Lemma.
From (\ref{eqdissection2}) and by Proposition \ref{propGamma/P} we deduce that for all $R>0$ it holds:
 \vspace{-3mm}        
 
\begin{equation}
\label{eqvdelta>}
 v^{4\Delta}_\Gamma (x,y,R) \geq   \frac12 \sum_{k= 0}^{N}   
  v^{\Delta/2}_{\Gamma} \!  \left(x,  k\Delta \right) \cdot v^{\Delta}_{P_i} \!  \left(x,y, (N-k)\Delta \right) 
\end{equation}
  
\noindent   as $\Delta  > 2 \ell_i$. Now, we set $h_i=b_{\xi_i} (x,y)$ and we sum  (\ref{eqvdelta>}) over annuli of radii $R_n = n \Delta$, and we get:
 \vspace{-5mm}       
       
\small  
$$
\hspace{-5mm} v_{\Gamma} (x,y, R) 
 \geq \frac14   \sum_{n=0}^{\lfloor \frac{R}{\Delta}\rfloor -2} 
               v_{\Gamma}^{4\Delta}  \left(x,y,  n\Delta \right)   
 \succ   \sum_{k=0}^{\lfloor \frac{R}{\Delta}\rfloor -1}  \left[ \sum_{n \geq k}^{\lfloor \frac{R}{\Delta}\rfloor -1} 
               v_{\Gamma}^{ \Delta/2} \left(x, (n-k)\Delta  \right) \right]    \cdot v_{P_i}^{\Delta} \left( x,y,  k\Delta \right)  \geq$$
\begin{equation}
\label{eqconvdiscr}
\hspace{3mm}    
 \geq   \sum_{k \geq \frac{h_i}{\Delta} +1}^{\lfloor \frac{R}{\Delta}\rfloor -1} \!\!\!\!\!\! 
                v_{\Gamma} \left(x,  R-(k+2) \Delta   \right)  \cdot v_{P_i}^{  \Delta} \left( x,y,  k\Delta \right) 
 \stackrel{C'}{\succ}  
           \! \sum_{k = \frac{h_i}{\Delta} +1 }^{ \lfloor \frac{R}{\Delta} \rfloor -1}  
               \!\!\! \frac{ v_{\Gamma} \left(x,  R- (k+2)\Delta  \right) }{ {\cal A}_{P_i} \left( x,   \frac{k\Delta + h_i}{2} \right)} 
 \end{equation}
 \normalsize

\noindent as   $v^{\Delta}_{P_i} (x,y,k\Delta) \succ {\cal A}^{-1}_{P_i} \left( x,   \frac{k\Delta + h_i}{2} \right)$
if $k\Delta \geq  h_i+\Delta>h_i+R_0$ by Proposition \ref{propVp}. 
 
\noindent Using again Proposition \ref{propVp} and (\ref{eqA(R+C)}), it is easily verified that the expression in (\ref{eqconvdiscr}) is greater than the continuous convolution 
 $v_{\Gamma}( x,  \cdot ) \ast v_{P_i}   ( x,y,  \cdot) \;(R+ 4\Delta)$, up to a multiplicative constant $CC' \Delta$.
This ends the proof, by taking $D_0 = 4\Delta$.$\Box$

 \vspace{5mm}
The Counting Formula enables us to reduce  the estimate of the growth function $v_X$ to a group-theoretical calculus, that is to the estimate of a the convolution of  $v_\Gamma$ with the cuspidal functions ${\cal F}_{P_i}$ of maximal parabolic subgroups  $P_i$  of $\Gamma$:

\begin{proposition}[Volume Formula]
\label{propvolume}
${}$

\noindent There exists  a constant $C'''=C'''(n,a,b,d,vol({\cal K}))$, such that:
\vspace{-3mm}

%

\small
\begin{equation}
\label{eqVx}
 \!\!  \left[ v_\Gamma (x, \cdot)  \ast  \!   \! \sum_i \!  {\cal F}_{P_i} (x,\cdot) \right]  \! \! ( R\!  - \!   2D_0  )    
\stackrel{C'''}{\prec}  
   v_X (x,R) 
\stackrel{C'''}{\prec}
   \left[ v_\Gamma (x, \cdot)  \! \ast \! \!  \sum_i \! {\cal F}_{P_i} (x,\cdot) \right]  \! \! \left( R \! \! +\!\!   2D_0 \right)
   \;\; \forall R\! \geq\!  0
\end{equation}
 \normalsize    


\noindent for $D_0=D_0(n,a,b,d)$ as in Proposition \ref{propdissection}.
\end{proposition}


\pagebreak

 \noindent {\bf Proof.} 
 Let $h_i = d(x, H_{\xi_i})$; we may assume that the constants $R_0, D_0$ of Propositions \ref{propVp} and \ref{propdissection} satisfy $D_0 \gg d \geq  \mbox{diam}({\cal K}) \geq h_i \gg R_0$.  
 Now call ${\cal S}_i (h) = \psi_{\xi_i, h} [{\cal S}_i]$; integrating $v_\Gamma (x,y,R)$ over the fundamental domain ${\cal D}$ yields,  by Proposition \ref{propdissection}:
\vspace{-5mm}

$$      \!\!\! v_X(x, R+2D_0)
  = \!\!  \int_{{\cal D}} \!\! v_{\Gamma} (x,y,R+2D_0) dy 
  =  \!\!   \int_{{\cal K}} \!\!  v_{\Gamma} (x,y, R+2D_0) dy \, + \sum_{i=1}^m \! \int_{{\cal C}_i} \!\! v_{\Gamma} (x,y,R+2D_0) dy$$

\vspace{-3mm}
     
$$   \stackrel{C''}{\succ}   
 \sum_{i=1}^m \int_{2h_i}^{R+D_0}
        \!\!\!\!\!\!\! v_\Gamma \left(x, R+2D_0-t \right)  
                     \left[   \int_{h_i}^{t-h_i} \!\!  \int_{{\cal S}_i(h)}   \!\! v_{P_i} \left(x,y, t  \right) dy\,   dh \right] dt$$
 


\noindent
which then  gives   by Propositions \ref{propVp} and  \ref{propFp},  as 
$ h= b_{\xi_i} (x,y)     \leq t - h_i   < t-R_0$,  
      $$
      \int_{2h_i}^{R+D_0} \!\!\!\!\!\!\!\!  v_\Gamma \left(x, R+2D_0 -t \right)  \left[   \sum_{i=1}^m 
                        \int_{h_i}^{t-h_i} \!\!  \frac{ {\cal A}_{P_i} (x,h)  }{{\cal A}_{P_i} \left(x, \frac{t+h}{2} \right)}  dh \right] dt$$

%
\vspace{-5mm}
   
$$ \geq  \int_{0}^{R+D_0-2h_i}  \!\!\!\!\!\!\!\!\!\!\!\!\!\!\!\!\!\!  v_\Gamma \left(x, R-t \right)   \left[    \sum_{i=1}^m 
                   \int_{0}^{t} \!\!  \frac{ {\cal A}_{P_i} (x,s+h_i)  }{{\cal A}_{P_i} \left(x, \frac{t+s+3h_i}{2} \right)}  ds \right] dt 
     \succ   \int_{0}^{R}  \!\!\!   v_\Gamma \left(x, R-t \right)   \sum_{i=1}^m  {\cal F}_{P_i}(x,t) dt       .$$    
                        

\noindent Reciprocally,  we have
  $v_{\Gamma} (x,R - D_0 ) \leq v_{\Gamma} (x,y,R)  \leq v_{\Gamma} (x,R+D_0)$
  so  again by  Proposition \ref{propdissection} and Remarks  \ref{reminequalities} we obtain
 \vspace{-5mm}
 
 %
 
 $$ v_X(x, R-2D_0) 
 \stackrel{C''}{\prec}   \mbox{vol} ({\cal K}) \cdot v_\Gamma (x, R-D_0)  
    +   \sum_{i=1}^m \int_{{\cal C}_i}  \left[  \int_0^{R-2D_0} \!\!\!\!\!\!\!\!\! v_\Gamma (x, t )    v_{P_i} (x,y,R-t) dt \right]  dy $$
 $$\hspace{17mm}  \stackrel{C'''}{\prec}  v_\Gamma (x, R-D_0)  
  +  \int_0^{R-2D_0} \!\!\!\!\!\!\!\!\!\! v_\Gamma (x,  t ) 
        \left[ \sum_{i=1}^m \int_{_0}^{R-t} \!\!\!  \frac{  {\cal A}_{P_i} (x,h)  }{   {\cal A}_{P_i} \left( x,\frac{R-t+h}{2} \right)  } dh \right] dt $$
as $v_{P_i} (x,y,R-t)=0$ for $R-t<b_{\xi_i}(x,y)=h$.  
This proves the converse inequality, 
since 
$v_\Gamma (x, R-D_0) \prec v_\Gamma (x, R-D_0) {\cal F}_{P_i} (R_0) \leq 
\frac{1}{D_0-R_0}\int_{R-D_0}^{R-R_0} v_\Gamma (x,t) {\cal F}_{P_i} (x,R-t)dt$.$\Box$

\vspace{5mm}        
As a consequence of the Volume Formula and of  Corollary \ref{cordeltaF}, we 
deduce\footnote{Part (i) of this  corollary already appears in  \cite{DPPS}, where an {\em upper} estimate  for $v_X$ is proved. \linebreak
Notice that in \cite{DPPS} we erroneously stated that also
$\omega^- (X) = \max\{ \delta(\Gamma), \omega^-( {\cal F}_{P_1}), ..., \omega^-({\cal F}_{P_m})\}$;
an explicit counterexample to this is given in Example \ref{exsparse}.}:

\begin{corollary}
\label{corvolume}
If  ${\cal F}_{P_i}$ are the cuspidal functions of the parabolic subgroups    of $\Gamma$:

\noindent (i) $\;\omega^+ (X) = \max\{ \delta(\Gamma), \omega^+( {\cal F}_{P_1}), ..., \omega^+({\cal F}_{P_m})\}$.

\noindent (ii)  $\omega^+ (X) = \omega^- (X) = \delta(\Gamma)$ if $\Gamma$ is $\frac12$-parabolically pinched.
 \end{corollary}

\section{Margulis function for regular lattices}
\label{margulis}

 In this section we assume that $\Gamma$ is a lattice which is neither sparse  nor exotic. \\
 To prove the  the divergence of the Poincar\'e series of $\Gamma$, 
   we  will need  a general criterion     which can be found in \cite{DOP}, \cite{DPPS2}:
  \vspace{3mm}
 
 \noindent {\sc Divergence Criterion}.  
 {\em 
Let $\Gamma$ be a  geometrically finite group: if $\delta^+(P) < \delta (\Gamma)$ \linebreak for every parabolic subgroup $P$ of $\Gamma$, then $\Gamma$ is divergent.
}
 
   \vspace{3mm}
 \noindent From the divergence, we will then deduce the finiteness of the Bowen-Margulis measure by the following result, due to Dal'Bo-Otal-Peign\'e (see \cite{DOP}):
  \vspace{3mm}

\noindent {\sc Finiteness Criterion}.  
{\em Let $\Gamma$ be a divergent, geometrically finite group,    $\bar X = \Gamma \backslash X$. We have $\mu_{BM} (U\bar X) <\infty$ if and only if for every maximal parabolic subgroup $P$ of $\Gamma$
\begin{equation}
\label{eqcritfiniteness}
\sum_{p\in P} d(x,px)e^{-\delta(\Gamma) d(x,px)} < +\infty.
\end{equation}
 }

 {\bf Proof of Theorem \ref{teor14basic}}.  
Let $\Gamma$ be a nonuniform lattice of $X$ which is neither sparse nor exotic. 
As $\Gamma$ is not exotic, it satisfies the  gap property  $\delta(P) < \delta(\Gamma)$ for all parabolic subgroups; by the Divergence and Finiteness Criterion recalled in \S1, we deduce that the group is divergent and that  $\mu_{BM} (U\bar X) <\infty$.
Therefore  $v_\Gamma (x,R) \stackrel{c_\Gamma (x)}{\asymp} e^{\delta (\Gamma) R}$ is purely exponential (for some $c_\Gamma (x)$ depending on $\Gamma, x$).  We will now show that $X$ has a Margulis function. \\ 
Let ${\cal D}$ be  the fundamental domain for $\Gamma$ and $P_i$   the maximal parabolic subgroup fixing $\xi_i$ as at the beginning of \S3: we call   $w (x,y,R) =  v_\Gamma (x,y,R) e^{-\delta(\Gamma)R}$, so that
  have

 \small
   \begin{equation}
\label{eqmargulis}
\frac{v_X (x,R)}{e^{\delta(\Gamma)R }}
=  \int_{{\cal D}}  \frac{v_\Gamma (x,y,R)}{e^{\delta(\Gamma)R}}dy 
    =  \int_{{\cal K}}  w(x,y,R) dy 
  +   \sum_{i=1}^m   \int_{{\cal C}_i} w(x,y, R)dy
\end{equation}
\normalsize       
  
\noindent We know that $v_\Gamma (x,y,R) \leq v_\Gamma (x,R+d)\leq c_\Gamma (x) e^{\delta(\Gamma)R}$ for $y \in  {\cal K}$, so we can pass to the limit for $R \rightarrow \infty$ under the integral sign in the first term. For the integrals over the cusps, we have: 
  \vspace{-4mm}       
  
  \small
  $$ w (x,y, R) 
  \stackrel{C''}{\prec}    \frac{ \left[ v_{\Gamma} (x, \cdot ) \ast v_{P_i} (x,y, \cdot) \right]   \!(R\!+\!D_0)}{e^{\delta(\Gamma)R} }  
  \stackrel{c_\Gamma (x)}{\prec}  \int_{  b_{\xi_i}(x,y)  }^{\infty}   \frac{ e^{-\delta (\Gamma) t}  }{ {\cal A}_{P_i} \left( x, \frac{ b_{\xi_i}(x,y) + t}{2}  \right) }dt   = w(x,y)
  $$ 
\normalsize


 \vspace{-2mm}        
 \noindent Notice that the dominating  function $w(x,y)$ is finite as $\delta^+(P_i) < \delta (\Gamma)$.\\
 We will now show that $w(x,y) \in L^1({\cal C}_i)$. With the same notations  $h_i= d(x, H_{\xi_i})$ and ${\cal S}_i (h) = \psi_{\xi_i, h} ({\cal S}_i)$ as before, we have for all $i$:
 \vspace{-5mm}         
  
 \small
$$\int_{{\cal C}_i}   w (x,y)  dy
  =\int_{h_i}^{\infty} \!\!\!  \int_{y \in {\cal S}_i (h) } \left[    \int_{  b_{\xi_i}(x,y)  }^{\infty}   \frac{ e^{-\delta (\Gamma) t}  }{ {\cal A}_{P_i} \left( x, \frac{ b_{\xi_i}(x,y) + t}{2}  \right) }dt   \right] dydh
   = \int_{h_i}^{\infty}  \!\!\!   \int_{h}^{\infty}   \frac{  e^{-\delta(\Gamma)t} {\cal A}_{P_i}  (h) }{ {\cal A}_{P_i}  \left( x, \frac{h+ t}{2}  \right) } dt dh
$$            

\vspace{-3mm}
                           
\begin{equation}
\label{eqdominatedconvergence2} 
\hspace{-10mm} 
 = \int_{h_i}^{\infty}  e^{-\delta(\Gamma) t }  \left[  \int_{h_i}^{t}   \frac{  {\cal A}_{P_i}  (h) }{ {\cal A}_{P_i}  \left( x, \frac{h+ t}{2}  \right) }dh \right] dt
 \stackrel{C}{\prec}  \int_{h_i}^{\infty}  e^{-\delta(\Gamma) t }  {\cal F}_{P_i} (t) dt  
\end{equation}      
\normalsize

\noindent which converges,  as $\Gamma$  is not sparse and so $\omega^+({\cal F}_{P_i})\leq \delta^+(P_i) < \delta (\Gamma)$,  by Corollary \ref{cordeltaF}.  
We therefore obtain from (\ref{eqmargulis}), by dominated convergence, using Roblin's asymptotics       
 $$\lim_{R \rightarrow +\infty} \frac{v_X (x,R)}{e^{\delta(\Gamma)R }} 
 =   \frac{\parallel \! \mu_x \! \parallel}{\delta(\Gamma) \parallel \! \mu_{BM} \! \parallel} \int_{{\cal D}} \!\!  \parallel \! \mu_y \!\parallel dy \;=:   m(x)< +\infty. $$

 \noindent  Notice that $m(x)$ defines an $L^1$-function on $\bar X = \Gamma \backslash X$, as its integral over ${\cal D}$ is finite.$\Box$
  \vspace{3mm}

{\bf Proof of Theorem \ref{teorex}(i)}.  
We assume now that  $X$ has an exotic lattice $\Gamma$, with the dominant parabolic subgroups $P_i$, for $i=1,...,d$, satisfying  $\delta:=\delta(\Gamma) =\delta^+ (P_i) \leq 2  (\delta^- (P_i) -\epsilon )$, for some $\epsilon>0$.
When  $\mu_{BM} (U\bar X)<\infty$,   the same lines of the above proof   apply: $v_\Gamma (x,R) \asymp c_\Gamma (x)e^{\delta R}$ is purely exponential, and for the same functions $w(x,y,R)$, $w(x,y)$   we again obtain  (\ref{eqdominatedconvergence2});
but we need some more work to deduce that, for the dominant cusps $P_i$,  the integral of $e^{-\delta t }  {\cal F}_{P_i} (t) $ converges.  
So, for every dominant subgroup $P_i$, we write
$v_{P_i} (x,t) = o_i (t)e^{\delta  t}$, for some subexpo\-nen\-tial  functions $o_i (t)$; so, 
${\cal A}_{P_i} (x,t) \asymp  e^{-2\delta t}/ o_i (2t)$  for $t\geq R_0$.
As $\Gamma$ is exotic, the   dominant parabolic subgroups $P_i$ are convergent: actually, for any divergent subgroup $\Gamma_0 < \Gamma$ with limit set $L(\Gamma_0) \subsetneq L(\Gamma)$ one has \linebreak 
$\delta(\Gamma_0 ) < \delta (\Gamma)$ (see \cite{DP}). Therefore, the Poincar\'e series of $P_i$ gives, for $\Delta > \Delta_0 \gg 0$  
$$      \infty > \sum_{p \in P_i} e^{-\delta d(x,px)} 
      \succ \sum_{k \geq 1}   \frac{v^{\Delta}_{P_i} (x,k\Delta)}{e^{\delta k}} 
       \asymp  \int_{\Delta}^{\infty} o_i(t) dt $$

 \noindent by Proposition \ref{propVp},  so the functions $o_i (t)$ are integrable.
 This shows that  
 $$w(x,y) = \int_{  b_{\xi_i}(x,y)  }^{\infty}   \frac{ e^{-\delta t}  }{ {\cal A}_{P_i} \left( x, \frac{ b_{\xi_i}(x,y) + t}{2}  \right) }dt 
  = e^{\delta b_{\xi_i}(x,y)}\int_{  b_{\xi_i}(x,y)  }^{\infty}    o_i (h+t) dt <\infty$$

 \noindent Moreover, as every dominant $P_i$ is  strictly $\frac12$-pinched, we  have 
$v_{P_i} (x,t) \succ   e^{\frac12 (\delta+ \epsilon) t}$ for some $\epsilon>0$, that is
${\cal A}_{P_i} (x,t) \prec  e^{-(\delta+\epsilon) t}$   for all $t>0$.
Then Proposition \ref{propFp} yields 
\vspace{-3mm}        
        
\begin{equation}
  \label{eqFo}
  {\cal F}_{P_i} (R) \asymp \int_{0}^R \frac{ {\cal A}_{P_i} (s) }{ {\cal A}_{P_i} (\frac{s+R}{2}) } ds
 \prec e^{\delta R} \int_{0}^R  e^{- \epsilon s}  o_i (s+R)  ds
 \hspace{1cm} \mbox{ for } R\gg 0
\end{equation}

\noindent hence   (\ref{eqdominatedconvergence2}) gives in this case:
\vspace{-3mm}        
        
\small
 $$  \int_{{\cal C}_i}   w (x,y)  dy
 \stackrel{C}{\prec}  \int_{h_i}^{\infty} \!\!\! e^{-\delta(\Gamma) t }  {\cal F}_{P_i} (t) dt  
 \asymp \int_{h_i}^\infty  \!\! \left[ \int_0^t  \!\!\! e^{-\epsilon s} o_i (s+t)ds  \right]  dt
 \leq \int_{0}^\infty  e^{-\epsilon s} \left[ \int_s^\infty  \!\!\! o_i (s+t)dt  \right]  ds
$$
\normalsize      
 
 \noindent which converges, since $o_i$ is integrable. We can therefore pass to the limit for $R\rightarrow \infty$  under the integral in  (\ref{eqmargulis}), obtaining the asymptotics for $v_X (x,R)$ as  before.

 \noindent On the other hand,  if $\mu_{BM} (U\bar X) = \infty$, then $v_\Gamma (x,R) = o_\Gamma (R)e^{\delta R}$ is lower-exponential, and by (\ref{eqFo}) we  have 
  ${\cal F}_{P_i} (x, R) = f_i (R) e^{\delta R}$ 
 with    $f_i (R)= \int_0^R  e^{- \epsilon s}  o_i (s+R)  ds$ for the dominant cusps, 
 and  $f_i (R) \prec  e^{- \epsilon R}  $, with $\epsilon >0$,  for the others; in both cases, 
 $f_i \in L^1$,  \linebreak since  the functions $o_i (t)$ are subexponential.  Proposition \ref{propvolume} then gives, for any arbitrarily small $\varepsilon' >0$

 \small
$$ \hspace{-10mm} \frac{v_X (x,R)}{e^{\delta R}}  
     \prec \frac{1}{e^{\delta R}}   \int_0^{R}  \!\! v_\Gamma (x,t) \sum_i {\cal F}_{P_i} (R -t) dt
    \prec  \int_0^{R} \!\!   o_\Gamma (t) \sum_i f_i (R -t) dt$$
    $$ \hspace{25mm} \leq  \sum_i   \parallel \! f_i \! \parallel_{_{1}} \cdot \sup_{t>\frac{R}{2}} o_\Gamma (t) \; +  \parallel \! o_\Gamma  \! \parallel_{\infty} \cdot \sum_i \int_{R/2}^{R}\!\! f_i (t)dt 
 \; \leq \; \varepsilon' \cdot \left( \sum_i \parallel \! f_i \! \parallel_{_{1}}  \!+\! \parallel \! o_\Gamma  \! \parallel_{\infty}  \right) $$
\normalsize

 \noindent  provided that $R\gg 0$, since  $o_\Gamma (t)$ is infinitesimal and  the $f_i$ are integrable. This shows that   $v_X(x,R)$  is lower-exponential too.$\Box$
 
 \begin{remark}
 {\em We have seen that, if $\mu_{BM} (U\bar X)=\infty$, then $v_\Gamma (x,R) = o_\Gamma (R)e^{\delta R}$ and $v_X (x,R) = o_X (R)e^{\delta R}$, where $ o_\Gamma, o_X$ are infinitesimal, and that ${\cal F}_{P_i} (x, R) \!\! = \!\!  f_i (R) e^{\delta R}$ with $f_i \in L^1$;  so,
 $$  \parallel \!  o_\Gamma   \! \parallel_{_{1}}   
 \prec \parallel \! o_X   \! \parallel_{_{1}} 
 \leq \int_0^\infty \frac{v_X (x,R)}{e^{\delta R}}dR  \prec \int_0^\infty  \int_0^R o_\Gamma(t) \sum_i f_i (R-t) dt dR
 \leq  \parallel \! o_\Gamma \! \parallel_{_{1}} \cdot \sum_i   \parallel \! f_i \! \parallel_{_{1}}   $$
           
\noindent and we can say that  $o_\Gamma $ is $L^1$ if and only if  $o_X $ is.
 }
 \end{remark}

Finally, in order to prove Theorem   \ref{teor14}, we  need to recall a characterization of constant curvature spaces as those  pinched, negatively curved spaces whose lattices realize the least possible value for the entropy. The minimal entropy problem  has a long history and has been declined in many different ways so far; see  \cite{knieper2}, \cite{BK},\cite{courtois} for the analogue of the following statement in the compact case, and \cite{dpps-rigidity} for a proof in the finite-volume case: 
\vspace{-2mm} 
 
\begin{theorem}
\label{teorbcgfinite}
Let  $\Gamma$  be  a  lattice  in a  Hadamard manifold $X$ with  pinched curvature $-b^2 \leq K_X\leq -a^2 <0$. Then  $\delta(\Gamma) \geq (n-1)a$, and $\delta(\Gamma)=(n-1)a$ if and only if  $X$ has constant curvature $-a^2$.
\end{theorem}

  {\bf Proof of Theorem \ref{teor14}.}
Assume that $\Gamma$ is a nonuniform lattice in a $\frac14$-pinched negatively curved manifold   $X$, i.e. $-b^2 \leq K_X \leq -a^2$ 
 with $ b^2  \leq 4 a^2$. If $X = \mathbb{H}^n_a$, then clearly $v_X (x,R) \asymp v_\Gamma (x,R)$ is purely exponential, $X$ has a Margulis function, and $\Gamma$ is divergent. 
Otherwise, let $P_i$ be the  maximal parabolic subgroups of $\Gamma$, up to conjugacy. By the formulas (\ref{eqA(R+C)}), we know that for all $x\in X$
$e^{-(n-1)bR}\prec {\cal A}_{P_i} (x, R)\prec e^{-(n-1)aR}$, 
so by Proposition \ref{propVp} we have
    \vspace{-3mm}    

$$\frac{a(n-1)}{2} \leq \delta^{-}(P_i) \leq \delta^{+}(P_i) \leq \frac{b(n-1)}{2}$$
for all $P_i$. 
Thus, $\Gamma$ is parabolically $\frac12$-pinched. It follows from Corollary \ref{corvolume}  that $\omega^+ (X) = \omega^- (X) = \delta (\Gamma)$. Moreover, for all $P_i$ we have
     \vspace{-3mm}   
 
$$\delta^+ (P_i) \leq  \frac{b(n-1)}{2} \leq a(n-1) < \omega (X) = \delta(\Gamma)$$

\noindent where the strict inequality follows by the rigidity Theorem \ref{teorbcgfinite}, since $X \neq \mathbb{H}^n_a$. \linebreak
The same argument applies when $\bar X$ is only asymptotically $\frac14$-pinched, by replacing $-a^2,-b^2$ with the bounds $-k_+^2  -\epsilon \leq K_X \leq -k_-^2  +\epsilon $ on the cusps $\bar {\cal C}_i$.
 Then,  $\Gamma$ is also non-exotic, and we can conclude  by Theorem \ref{teor14basic} that
 $\Gamma$ is divergent, with  finite Bowen-Margulis measure,  $v_X \asymp v_\Gamma$ and $X$ has a $L^1$ Margulis function $m(x)$.$\Box$

\section{Examples}
\label{examples}
In this section we show that all the cases presented in  Theorem \ref{teorex} do  occurr,   by providing examples of spaces $X$ with exotic or sparse lattices  $\Gamma$ which do not admit a Margulis function, and  with functions $v_{\Gamma}, v_X$ having different behaviour.

\noindent If $\bar {\cal  C}= P \backslash H_\xi (o)$ is a cusp of $\bar X \!=\! \Gamma \backslash X$, we write the metric  of $X$ in   horospherical co\-or\-di\-nates on \nolinebreak $H_{\xi} (o) \! \cong \!\partial H_{\xi} (o) \! \times \!\mathbb{R}^+$  as 
$g = T(x,t)^2 dx^2 + dt^2$, for $x \!\in\! \partial H_{\xi} (o)$ and \nolinebreak $t \!=\! b_{\xi} (o,\cdot)$. \linebreak We call the function  $T(x,t)$ the {\em analytic profile of the cusp} $\bar {\cal  C}$. The horospherical area ${\cal A} _P (x,t)$ is then obtained by integrating $T^{n-1}(x,t)$  over a compact fundamental  domain ${\cal S}$ for  the action of $P$ on  $\partial H_{\xi} (o)$; thus, we have  
\vspace{-3mm}

$$ \hspace{0mm}  {\cal A}_{P} (x,t) \stackrel{c}{\asymp} T^{n-1} (x, t)  \hspace{15mm} \mbox{for all }x \in \bar {\cal  C}$$

\noindent (for a constant $c$ depending on $X$ and $o$). Also, notice that, in the particular case where $T(y,t)=T(t)$, for points $x,y$ belonging to a same horosphere $H_\xi$ we have by the Approximation Lemma \ref{lemmapp} 
\vspace{-3mm}        
        
\begin{equation} 
\label{eqdistanceprofile}
\hspace{10mm}
d(x,y)  \sim  2 T^{-1} \left( \frac{T(0)}{d_\xi (x,y) } \right)  
\hspace{10mm} \mbox{for }  R =d(x,y) \rightarrow \infty.
\end{equation}

\vspace{3mm}
We will  repeatedly make use of the following lemma, which is a easy modification of one proved in \cite{DPPS}:
 \vspace{-3mm}       

\begin{lemma}
\label{lemmainterpolation}
Let $b>a>0$, $\beta>\alpha>0$ and $ \epsilon>0$   be given. \\ There exist $D=D(a,b, \alpha, \beta, \epsilon)>1$ and $D'=D'(a,b, \alpha, \beta)>0$ such that if   $[p,q]$, $[r,s]$ are disjoint  intervals satisfying   $r\geq Dq$ and $p\geq D'$, then  there exist  $C^2$, convex and decreasing functions $\phi_\epsilon$, $\psi_\epsilon$ on $[p,s]$ satisfying:
\vspace{-5mm}                

$$  \left\{\begin{array}{clll}
               \forall\, t\in[p,q],  & \phi_\epsilon (t)  =t^\beta e^{-b t}  										  \\
               \forall\, t\in[r,s],  & \phi_\epsilon (t)  =t^\alpha e^{-a t}  										 \\
               \forall\, t\in[p,s],  & t^\beta e^{-b t}  \leq  \phi_\epsilon (t)  \leq t^\alpha e^{-a t}  							 \\
             	 \forall\, t\in[p,s],    &  a^2-\epsilon\leq \frac{ \phi_\epsilon ''(t) }{ \phi_\epsilon (t) } \leq b^2+\epsilon    
         \end{array}\right. \mbox{ and }\;\;
      \left\{\begin{array}{clll}
               \forall\, t\in[p,q],  & \psi_\epsilon (t)  =t^\alpha  e^{-a t}  										  \\
               \forall\, t\in[r,s],  & \psi_\epsilon (t)  =t^\beta  e^{-b t}  										 \\
               \forall\, t\in[p,s],  & t^\beta e^{-b t}  \leq  \psi_\epsilon (t)  \leq t^\alpha e^{-a t}  							 \\
             	 \forall\, t\in[p,s],    &  a^2-\epsilon\leq \frac{ \psi_\epsilon ''(t) }{ \psi_\epsilon (t) } \leq b^2+\epsilon    
         \end{array}\right. $$
     
\end{lemma}

\vspace{2mm}
\begin{example}Sparse lattices. 
\label{exsparse}
{\em

\noindent Sparse lattices satisfying $\omega^+ (X) > \delta(\Gamma)$ were constructed by the authors in \nolinebreak \cite{DPPS}.  
Here, we modify that construction to show that, for spaces $X$ admitting sparse lattices, one can  have $\omega^+ (X) > \omega^- (X) > \delta(\Gamma)$  (in contrast, notice that $\delta(\Gamma)$   always is a true limit); this shows in particular  that sparse lattices generally do not have a Margulis function.\\
We start  from a hyperbolic surface
$\bar X_0 = X_0\backslash \Gamma$ of finite volume,  homeomorphic to a 3-punctured sphere, 
and, for any arbitrary small $\epsilon>0$,  we perturb  the hyperbolic metric $g_0$ on  one cusp $\bar{ {\cal C}}=P \backslash H_\xi (x)$  into a metric $g_\epsilon$  by choosing an analytic profile  $T_{\epsilon}$ obscillating, on infinitely many   horospherical bands, from $e^{-t}$ to $e^{-bt}$. 

\noindent Namely,  choose $a=1, b>2$ and $\epsilon>0$ arbitrarily small, and let $D, D'$ be the constants given by Lemma  \ref{lemmainterpolation}. 
 For $M \gg 1$, we   define a sequence of disjoint  subintervals of $[M^{4n}, M^{4n+1}]$:  
  \vspace{-3mm}
         
 \small
$$[p_{n},q_{n}] :=  [ M^{4n},  2 M^{4n}],  \hspace{1cm} [r_{n}, s_{n}] :=  \left[ \frac{p_n+M^{4n+1}}{2},\frac{q_n+M^{4n+1}}{2}\right]$$
\normalsize

\noindent such that $r_n \geq Dq_n$, $p_{n+1} \geq Ds_n$, \nolinebreak$p_1 \geq D'$ 
(we can choose  any $M \geq \max\{4D-1, \sqrt[3]{D}\}$   in order  that these conditions are satisfied). Notice that $\frac{t+M^{4n+1}}{2} \in [r_n,s_n]$ for all $t\in [p_n, q_n]$.
Then, by Lemma \ref{lemmainterpolation}, we consider a  $C^2$, decreasing function $T_\epsilon (t)$  satisfying:

\noindent (i)  $T_{\epsilon} (t) = e^{-t}$ for $t \in  [M^{4n-2},M^{4n}] \cup [p_n,q_n]$, and  $T_{\epsilon} (t) = e^{-bt}$ for $t \in [r_{n},s_{n}]$;

\noindent (ii) $e^{-bt} \leq T_{\epsilon} (t) \leq e^{-t}$ and 
 $-b^2-\epsilon \leq T_\epsilon '' (t) / T_\epsilon (t)  \leq -1+\epsilon$.  

\noindent  Thus, the new analytic profile $T_{\epsilon} (t)$ of $\bar {\cal C}$  coincides with the profile  of a usual hyperbolic cusp  on $[M^{4n-2}, 2M^{4n}]$, and with the profile  of a cusp in curvature $-b^2$ on the bands  $[r_{n},s_{n}] \subset  [M^{4n}, M^{4n+1}]$.
 We have, with respect to  the metric $g_\epsilon$:
 \vspace{2mm}        
 
\noindent {\bf (a)} $\delta^+(P) = \frac{b}{2}$ and  $\delta^-(P)=\frac12$,   by (i) and (ii),  because of Proposition \ref{propVp};
 \vspace{2mm}
 
\noindent {\bf (b)} $\omega^+ ({\cal F}_{P}) \geq \frac{b}{2} + \delta$ for  $\delta=\frac{1}{M}(\frac{b}{2}-1)>0$, because
for $R=M^{4n+1}$
 \vspace{-5mm}

\begin{equation}
\label{eqFmagg}
{\cal F}_{P} (x, R) 
\succ \int_{0}^{R} \frac{{\cal A}_\epsilon (x,t)}{{\cal A}_\epsilon (x, \frac{t+R}{2})}dt
\geq \int_{p_{n}}^{q_{n}} \frac{e^{-t}}{ e^{-b (\frac{t+R}{2}) }}dt 
\succ  e^{\frac{b}{2}R} \cdot M^{4n}e^{(\frac{b}{2}-1)p_{n}} 
\geq e^{\frac{b}{2}R}\cdot e^{ \frac{1}{M}(\frac{b}{2}-1)R  } 
\end{equation}

\vspace{-3mm}       
\noindent as $p_{n} /  R =\frac{1}{M}$;
         
\vspace{2mm}  
\noindent {\bf (c)} $\omega^- ({\cal F}_{P}) \leq \frac12$ if $M >2$, as  for $R \in [M^{4n+3}, M^{4n+4}]$ we obtain:
\vspace{-3mm}

\begin{equation}
\label{eqFmin}
{\cal F}_{P} (x,R) 
 \prec \int_0^R \!\!\!\!  \frac{  e^{-t} }{ e^{- (\frac{t+R}{2})} }  dt  
 \prec   e^{\frac{R}{2}}  
\end{equation}

\vspace{-2mm}
\noindent  since $M^{4n+4} \geq\frac{t+R}{2} \geq \frac{M^{4n+3}}{2} \geq M^{4n+2}$;

\vspace{2mm}
 \noindent {\bf (d)} $\delta(\Gamma)$ is arbitrarily close to $\delta^+(P)$,  let's say $\delta(\Gamma)\leq\frac{b}{2}+\frac{\delta}{2}$,  if we perturb  the hyperbolic metric sufficiently far in the cusp $\bar{ {\cal C}}$, i.e. if $r_1\gg0$  (this is  Proposition 5.1 in \cite{DPPS}).
 
 \vspace{2mm}
\noindent  It follows that $\omega^- (X) > \delta(\Gamma)$. Actually, assume that $v_\Gamma (x,R) \succ e^{(\delta(\Gamma)-\eta)R}$, for  arbitrarily small $\eta$. 
By Proposition \ref{propvolume}  and (\ref{eqFmagg}),  we deduce that 
for any $R\gg0$, if  $M^{4n+1} \leq R  < M^{4n+5}$
\vspace{-5mm}        
  
$$v_X (x,R+2 D_0) 
\geq   v_{\Gamma}  (x, \cdot)  \ast {\cal F}_{P} (x,  \cdot ) \; (x,R) 
\succ  e^{ (\delta(\Gamma)-\eta) (R- M^{4n+1})} \cdot e^{(\frac{b}{2} + \delta)M^{4n+1}}$$

  \noindent by  taking just  the  term $v_{\Gamma}  (x, R-t)   {\cal F}_{P} (x,  t))$ of the convolution with $t$ closest to $M^{4n+1}$, where ${\cal F}_{P} (t)  \succ e^{(\frac{b}{2}+\delta)t}$; as  $M^{4n+1} \geq R/M^4$ we get 
$v_X (x,R+2\Delta)  \succ e^{(\delta(\Gamma)  - \eta + \frac{\delta/2 +\eta}{M^4})R} $
 which gives
$\omega^- (X) \geq \delta(\Gamma)+ \frac{\delta}{2M^4}$, $\eta$ being arbitrary.
        
  \vspace{2mm}
\noindent Finally, we show that  $\omega^+ (X) >   \omega^- (X)$. In fact, the  cusps different from $\bar{ {\cal C}}$ being hyperbolic,  we have, always by Proposition \ref{propvolume},  that $\omega^+ (X)=\omega^+ ({\cal F}_{P}) \geq \frac{b}{2}+\delta$.

\noindent  On the other hand, we know that $\omega^+ ({\cal F}_{P}) \leq \max\{ \delta^+(P),  2 (\delta^+(P) -   \delta^-(P)\} = b-1$, by Corollary \ref{cordeltaF}; thus,  
assuming ${\cal F}_P (x,t) \prec e^{(b-1+\eta)t}$, for  arbitrarily small $\eta$, equation 
 (\ref{eqFmin}) yields for $R=M_{4n+4}$ 
$$     v_X (x,R-2D_0) 
\leq \int_{0}^{M^{4n+3}}  \!\! \!\!\! \!\!\!\! v_\Gamma (x,R-t)  \cdot {\cal F}_{P}(x,t) dt 
       + \int_{  M^{4n+3}}^{R} \!\! v_\Gamma (x,R-t)  \cdot  {\cal F}_{P}(x,t)dt$$ 

$$\hspace{10mm} \prec  \int_{0}^{M^{4n+3}}   \!\!\!\!\!\!  e^{\delta(\Gamma)(R-t)} \cdot e^{(b-1+\eta)t} dt
        +   \int_{  M^{4n+3}}^{R} \!\! e^{\delta(\Gamma)(R-t)}\cdot  e^{\frac12 t } dt$$


$$ \prec e^{\delta(\Gamma)R} \cdot e^{ (b-1 +\eta- \delta(\Gamma))   M^{4n+3}} \;\;  \leq \;\; e^{(\frac{b}{2}+ \frac{\delta}{2} + \frac{(b/2+\eta -1}{M} ) R}$$


\noindent being $\frac{b}{2} \leq \delta(\Gamma) \leq \frac{b}{2} + \frac{\delta}{2}$ and $ M^{4n+3} =  \frac{R}{M}$.
Hence  $\omega^- (X) < \frac{b}{2} + \delta \leq \omega^+ (X)$, if $M \gg 0$ and $\eta$ small enough.

}
\end{example}

\begin{examples} Exotic, strictly  $\frac12$-parabolically pinched  lattices.
\label{exexotic<}
{\em

\noindent We say that a lattice $\Gamma$ is  {\em strictly   $\frac12$-parabolically pinched} when every parabolic sugroup $P < \Gamma$ satisfies the strict inequality $\delta^{+}(P) < 2\delta^{-} (P)$.
Let $\bar X = \Gamma   \backslash X$ as before; we  show here that, for $\Gamma$ exotic and     strictly   $\frac12$-parabolically pinched,   the following cases which appear in Theorem \ref{teorex} do occur:

\noindent {\bf (a)} $\mu_{BM} (U\bar X) = \infty $ and $v_X$ is lower-exponential;

\noindent {\bf (b)} $\mu_{BM} (U\bar X) < \infty $ and  $v_X$ is purely exponential.


\vspace{2mm}
 We start by an example of lattice satisfying {\bf (a)}. \\
 In \cite{DPPSnew} the authors show how to construct {\em convergent  lattices},  in pinched negative curvature and any dimension $n$; we will take $n=2$ here by the sake of simplicity. \linebreak
 In those examples, the metric   is hyperbolic everywhere   but one cusp ${\cal C}$,  which has analytic profile
 $T(t) = t^\beta e^{bt}$ for $t\geq t_0 \gg 0$, with $\beta>1$ and $b>2$.
Therefore,   there is just one dominant  maximal parabolic subgroup $P$, with  ${\cal A}_{P} (x,t) \asymp T (t) \asymp e^{bt}$, and $\delta^+(P)=\delta^-(P)=\frac{b}{2}$; moreover, the subgroup $P$ is convergent as
\vspace{-5mm}

\small
$$     \sum_{p \in P} e^{-\frac{b}{2} d(x,px)} 
 \leq \sum_{k \geq 0}  v_P (x, k) e^{-\frac{b}{2} k} 
 \asymp  \int_1^\infty \!\!\! \!  \frac{  e^{- \frac{b}{2} t}    }{ {\cal A}_{P} (x,\frac{t}{2}) } dt
 \asymp  \int_1^\infty \!\!\! \!  \frac{  e^{- \frac{b}{2} t}    }{  (t)^\beta \cdot e^{-b \frac{t}{2}}    } dt 
\asymp \int_1^\infty \!\!\!\!\!\! t^{-\beta}dt < \infty.$$
 \normalsize    

\noindent By decomposing the elements of $\Gamma$ in geodesic segments which, alternatively,  either go very deep in the cusp or stay  in the hyperbolic part of $X$, we show in \cite{DPPSnew} that $\Gamma$ is convergent too, provided that $t_0 \gg 0$. 
Then, $\Gamma$ is exotic with  infinite Bowen-Margulis measure, and  $v_\Gamma (x,R)$ is lower-exponential by Roblin's asymptotics. By  Theorem \ref{teorex}(i),  the  function $v_X$  is  lower-exponential  as well,  with the same exponential growth rate.

 \vspace{2mm}
  We now give  an example for {\bf (b)}.  \\
  This is more subtle, as we need to take a {\em divergent, exotic lattice} $\Gamma$: the existence of such lattices is established, in dimension 2, in \cite{DPPSnew}. 
Again, the simplest example is homeomorphic to a three-punctured sphere, with  three cusps, and hyperbolic metric outside one cusp $\bar {\cal C}$, which has analytic profile 
 \vspace{-3mm}
 
 $$ 
T(t) =\left\{
\begin{array}{ll}
  e^{-t}& \mbox{ for }  t\leq A     \\
   e^{-b t}      &    \mbox{ for }  t \in [A, A+  B] +D      \\
 t^3  \cdot e^{-b t}  &   \mbox{ for } t \gg    D+ A + B 
\end{array}
\right.
$$
with   $b>2$ and  $A, B, D \gg 0$. 
As before,  we have one dominant and convergent maximal parabolic subgroup $P$, with  $\delta^+(P)=\delta^{-}(P)=  \frac{b}{2}$.  In \cite{DPPSnew} it is proved that, according to the values of $A$ and $B$, the behaviour of the group $\Gamma$ is very different: it is  convergent with critical exponent  $\delta(\Gamma) =  \delta^+(P)$, for $A\gg0$ and $B=0$, while it is divergent with $\delta(\Gamma) >  \delta^+(P)$ if  $B \gg A$. By perturbation theory of transfer operators, it is then proved that there exists a value of $B$ for which  $\Gamma$ is divergent with $\delta(\Gamma) =\delta^+(P)$ precisely. 
Thus, for this particular value of $B$, the lattice $\Gamma$  is exotic, and has finite Bowen-Margulis measure by the Finiteness Criterion, as   
\vspace{-5mm}

\small
\begin{equation}
\label{eqBM}
\sum_{p \in P} d(x,px)e^{-\delta(\Gamma) d(x,px)}     
\prec   \int_{1}^\infty \!\!\! \frac{  t e^{- \frac{b}{2}  t}    }{ {\cal A}_{P} (x,\frac{t}{2}) } dt
\prec   \int_{1}^\infty\!\!\! \frac{  t e^{- \frac{b}{2}  t}    }{  t^{3}  \cdot e^{-b\frac{t}{2}}} dt
\asymp   \int_{1}^\infty \!\!\! t^{-2} dt < \infty 
\end{equation}
\normalsize

 }
\end{examples}

\noindent It follows that $v_X \asymp v_\Gamma$ is purely exponential, by Theorem \ref{teorex}(i).

\begin{examples} Exotic, exactly $\frac12$-parabolically pinched  lattices.
\label{exexotic=}
{\em

\noindent We say that a lattice $\Gamma$ is  {\em exactly   $\frac12$-parabolically pinched} when it is $\frac12$-parabolically pinched and has a parabolic sugroup $P < \Gamma$ satisfisfying the quality $\delta^{+}(P) = 2\delta^{-} (P)$.
We show here that  for an exotic and exactly $\frac12$-parabolically pinched lattice $\Gamma$,   the following cases can occur:

\noindent {\bf (a)} $\mu_{BM} (U\bar X) < \infty $,   with $v_\Gamma$ purely exponential and $v_X$ upper-exponential;

\noindent {\bf (b)} $\mu_{BM} (U\bar X) = \infty $, with $v_\Gamma$ lower-exponential and $v_X$ 
upper-exponential.

 
 \vspace{3mm}
  We start by  {\bf (a)}. Consider a surface with three cusps  as in the Examples \ref{exexotic<}, now perturbing the hyperbolic metric   on the cusp $\bar {\cal C}$ to an analytic profile defined as follows.
First,  choose  a sequence of disjoint subintervals of $[M^{2n}, M^{2n+1}]$
  \vspace{-5mm}      
        
\begin{equation}
\label{eqintervals}
[p_{n},q_{n}] :=  [ M^{2n},  \mu M^{2n+1}],  
\hspace{1cm} 
   [r_{n}, s_{n}] :=  \left[ \frac{p_n+M^{2n+1}/2}{2},\frac{q_n+M^{2n+1}}{2} \right]
\end{equation}
      
 \normalsize  
\noindent and then define, for $b >1$ and $0< \gamma <1$      
 $$ 
T (t) =\left\{
\begin{array}{ll}
  e^{-t}& \mbox{ for }  t\leq A     \\
   e^{-b  t}      &    \mbox{ for }  t \in  [A, A+B]  +D     \\
   t   \cdot e^{-\frac{b }{2} t}  &   \mbox{ for } t \in  [p_n,q_n]  \\
   t^{2+\gamma} \cdot e^{-b  t}  &   \mbox{ for } t \in  [r_n,s_n]  \\
\end{array}
\right.
$$
 \noindent 
with  $t^{2+\gamma}e^{-bt} \leq T (t) \leq t   \cdot e^{-\frac{b }{2} t}  $ for all $t \geq t_0 \gg 0$
(in order that the conditions of Lemma \ref{lemmainterpolation} are satisfied, it is enough to  choose 
  any $0<\mu < \frac{1}{4D}$  and  $M > D$).  
 
 \noindent   As before, the profile $T$ gives a divergent, exotic lattice $\Gamma$ for a suitable value of $B$ and $A\gg0$, with dominant parabolic subgroup $P$ having $\delta^+(P) =\frac{b}{2}= \delta(\Gamma)$, and $\delta^-(P) =\frac{b}{4}$.
  The Bowen-Margulis measure of $\Gamma$ is finite, as (\ref{eqBM}) also holds  in this case; thus, $v_\Gamma$ is purely exponential. Let us now show that $v_X$ is upper exponential: for every $R=M^{2n+1}$ we have, by 
 Proposition \ref{propvolume}, 
 $$ v_X(x, R +2D_0) 
 \succ  \left[ v_X (x, \cdot) \ast {\cal F}_{P} (x, \cdot)\right]  (R)
 \asymp \int_0^R v_\Gamma (x, R-t) 
           \left[     \int_0^t \frac{ {\cal A}_{P} (x,s)     }{    {\cal A}_{P} (x,\frac{s+t}{2})     } ds \right] dt$$
 $$= \int_0^R {\cal A}_{P} (x,s)  
      \left[   \int_s^R  \frac{ v_\Gamma (x, R-t)   }{    {\cal A}_{P} (x,\frac{s+t}{2})     }  dt  \right] ds 
\geq\int_{p_n}^{q_n}   {\cal A}_{P} (x,s)  
      \left[   \int_{\frac{R}{2} }^{R}  \frac{ v_\Gamma (x, R-t)   }{    {\cal A}_{P} (x,\frac{s+t}{2})     }  dt  \right] ds$$

\noindent since $q_n <\frac{R}{2}$. As $\frac{s+t}{2} \in [r_n,s_n]$ if $s\in [p_n,q_n]$ and $t\in [\frac{R}{2} ,R]$, by  the definition of $T (t) \asymp {\cal A}_{P}(x, t)$ on $[r_n,s_n]$, this yields
\vspace{-3mm}         
 
 \small
 $$ v_X(x, R)
  \succ \int_{p_n}^{q_n}   s e^{-\frac{b}{2}s}   
      \left[   \int_{\frac{R}{2}}^{R}  
      \frac{  e^{\frac{b}{2}(R-t)}   }{  e^{- b ( \frac{s+t}{2})} (s+t)^{2+\gamma}   }  dt  \right] ds
 \succ e^{\frac{b}{2} R}  \int_{p_n}^{q_n} \frac{Rs}{(s+R)^{2+\gamma}} ds
      $$
     
\normalsize
\noindent with 
$\displaystyle  \int_{p_n}^{q_n} \frac{Rs}{(s+R)^{2+\gamma}} ds 
\geq  \int_{\frac{1}{M}}^{\mu} \frac{  u }{   (1+u)^{2+\gamma}}du \asymp R^{1-\gamma}$, so $v_X$ is upper-exponential.

\pagebreak
Producing examples for case {\bf (b)}  is more difficult; for this, we will need an exotic lattice $\Gamma$ whose  orbital function satisfies $v_\Gamma (o,R) \asymp \frac{1}{R^\gamma}e^{\delta(\Gamma)R}$. 
Lattices with lower-exponential growth and infinite Bowen-Margulis measure are investigated in \cite{DPPSnew}, where a refined counting result is proved, according to the behaviour of the profile functions of the cusps (the examples in \cite{DPPSnew} are, as far as we know,  the only precise estimates of the orbital function for groups with infinite Bowen-Margulis measure). Here we  only give the necessary analytic profiles of the cusps  in order to have a function $v_X$ which is exponential or upper-exponential,  referring to \cite{DPPSnew} for the  precise estimate of $v_\Gamma$.

\noindent  We again start from a hyperbolic surface $\bar X_0 = X_0\backslash \Gamma$  with three cusps as in   \ref{exexotic<}, and  perturb now the metric on {\em two} cusps. We choose $b>2$ and  $ 1+\gamma <\beta <2+\gamma $, and define the profiles for $\bar {\cal C}_1$ and $\bar {\cal C}_2$ as
\vspace{-5mm}  
        
 $$ 
 T_1 (t) =\left\{
\begin{array}{ll}
  e^{-t}& \mbox{ for }  t\leq A     \\
   e^{-b  t}      &    \mbox{ for }  t \in  [A, A+B]  +D     \\
   t    \cdot e^{-\frac{b }{2} t}  &   \mbox{ for } t \in  [p_n,q_n]  \\
   t^\beta \cdot  e^{-b  t}  &   \mbox{ for } t \in  [r_n,s_n]  \\
\end{array}
\right. 
\; \mbox{ and } \;\;\;
T_2 (t) =\left\{
\begin{array}{ll}
  e^{-t}& \mbox{ for }  t\leq A     \\
   t^{1+\gamma} e^{-bt}    &    \mbox{ for }  t \gg A     \\
  \end{array}
\right. 
$$

 \noindent for the same sequence of intervals $[p_n, q_n]$, $[r_n, s_n]$ as in (\ref{eqintervals}). 

\noindent  If $P_1, P_2$ are the   associated maximal parabolic subgroups,  we have  $\delta^- (P_1) =  \frac{b}{4}$ and $\delta^+ (P_1) =  \frac{b}{2} $, while  $\delta^+ (P_2) = \delta^- (P_2)= \frac{b}{2}$ by construction. It is easily verified that these parabolic subgroups are convergent as $\gamma>0$.
Again, pushing the perturbation far in the cusps (i.e. choosing $A \gg 0$) and for a suitable value of $B$, the lattice $\Gamma$ becomes exotic and divergent; it has two dominant cusps,  it is exactly $\frac12$-parabolically pinched,  and has  infinite Bowen-Margulis measure, because (as $\gamma<1$)
\vspace{-3mm}   

$$\sum_{p \in P_2} d(x,px)e^{-\delta(\Gamma) d(x,px)}     
\prec   \int_{1}^\infty\!\!\! \frac{  t e^{- \frac{b}{2}  t}    }{  t^{1+\gamma}  \cdot e^{-b\frac{t}{2}}} dt
\asymp   \int_{1}^\infty \!\!\! t^{-\gamma} dt = \infty.
$$

\noindent Accordingly, $v_\Gamma$ is lower-exponential. 
 In \cite{DPPSnew} it is proved that the {\em least convergent} dominant parabolic subgroup determines the asymptotics of $v_\Gamma$; in this case, the parabolic subgroup $P_1$ converges faster than $P_2$, and the chosen profile for $\bar {\cal C}_2$ then gives $v_\Gamma (o,R) \asymp \frac{1}{R^{1-\gamma}}e^{\delta(\Gamma)R}$,  {\em provided that $\gamma \in (\frac12, 1)$}, cp. \cite{DPPSnew}.

%
%
%
%
%
 
\noindent  Let us now estimate $v_X (x,R)$, for $R=M^{2n+1}$. 
 Writing 
 $T_1 (t) = \tau^+(t) e^{-bt}=  \tau^-(t) e^{-\frac{b}{2} t}$ 
 so that $\tau^{+}(t)=t^\beta$ on $[r_n, s_n]$ and  $\tau^{-}(t)=t$ on $[p_n, q_n]$, we  compute as in case \textbf{(a)}:
 \vspace{-1mm}       
  
  \small         
$$\hspace{-18mm} v_X( x, R+2D_0) 
\succ \left( v_\Gamma (x, \cdot) \ast {\cal F}_{P_1} (x, \cdot)  \right) (R)
= \int_0^R \int_0^t   
    \frac{  {\cal A}_{P_1} (x, s)  }{   {\cal A}_{P_1} (x, \frac{t+s}{2} ) } 
     v_\Gamma (x,R-t) dtds$$
$$\hspace{14mm}  \asymp   \int_0^R \!\!\! \int_0^t  
 \frac{   \tau^-(s) \cdot e^{-\frac{b}{2}s} \cdot e^{\frac{b}{2}(R-t)}    }
      {  \tau^+(\frac{t+s}{2})\cdot (R-t)^{1-\gamma} \cdot e^{-b(\frac{t+s}{2})}   } dtds   
=   e^{\frac{b}{2}R} \int_0^R  \!\!\!  \tau^-(s) 
 \left[ \int_s^R \!\!\! \frac{ dt }{ \tau^+(\frac{t+s}{2}) (R-t)^{1-\gamma} }  \right]  ds  $$ 
$$\hspace{-15mm} \succ   e^{\frac{b}{2}R} \int_{p_n = \frac{R}{M} }^{q_n =\mu R} s  \left[  \int_{\frac{R}{2}}^R  \frac{ dt  }{  R^\beta  ( R-t )^{1-\gamma} }   \right]  ds  \succ \left(\mu - \frac{1}{M} \right) R^{2+\gamma-\beta} e^{\frac{b}{2}R}   $$ 
 \normalsize     
 
 \noindent which is upper-exponential as $\beta<2+\gamma$. 

}
\end{examples}

 \begin{remark}
 {\em Notice that in all these examples $b$ can be chosen arbitrarily close to $2a=2$. Thus, by the last condition in Lemma 
\ref{lemmainterpolation}, the analytic profiles give metrics with curvature
$-4a^2 - \epsilon \leq K_X \leq -a^2$, for arbitrarily small $\epsilon >0$. 
}
 \end{remark}



\end{document}